%% file: conference_101719.tex
\documentclass[conference]{IEEEtran}
\IEEEoverridecommandlockouts
\usepackage{cite}
\usepackage{dsfont}
\usepackage{amsmath,amssymb,amsfonts, comment,tikz}
\usepackage{mathptmx, multirow}
\usepackage{amsthm}
\usepackage{epsfig} 
\usepackage{bm}
\usepackage{slashbox}
\usepackage{algorithmic}
\usepackage{graphicx, nicefrac}
\usepackage{textcomp}
\usepackage{xcolor}
\usepackage{hyperref}
\usepackage{enumerate}

\newcommand{\ignore}[1]{}

\usepackage{fancyhdr}
\usepackage{lastpage}
\def\BibTeX{{\rm B\kern-.05em{\sc i\kern-.025em b}\kern-.08em
    T\kern-.1667em\lower.7ex\hbox{E}\kern-.125emX}}
    
\newtheorem{thm}{Theorem}

\newtheorem{defn}{Definition}
\newtheorem{lemma}{Lemma}
\newcommand{\eop} {\hfill{$\blacksquare$}}
\newcommand{\Prob}{{\mathbb{P}}}
\newcommand{\E}{{\mathbb{E}}}
\newcommand{\eL}{{\phi_L^*}}

\newcommand{\eU}{{\phi_U^*}}

\newcommand{\phifp}{\phi^*_{P}} 
\newcommand{\phimr}{\eU} 
\newcommand{\phiunderbar}{\underline{ \phi}}
\newcommand{\phibar}{{\bar  \phi}} 
\newcommand{\phiistar}{{  \phi_i^*}} 
\newcommand{\phido}{{\phi^o_d}}

\newcommand{\VPhi}{{\bm  \phi}}
\newcommand{\PB}{{\cal B}}

\newcommand{\DA}{{\cal D}}
\newcommand{\DAl}{{\tilde{\cal D}}}
\newcommand{\MR}{{\cal M}}
\newcommand{\eiw}{{e}}

\newcommand{\eow}{{e}}
\newcommand{\etw}{{e}}
\newcommand{\niw}{{n_{i}}}
\newcommand{\nir}{{r_{i}}} 

\newcommand{\mstar}{{m^*}}
\newcommand{\MRl}{{\tilde{\cal  M}}}
\newcommand{\PBl}{{\tilde{\cal B}}}
\renewcommand{\P}{{\cal P}}

\newcommand{\lamol}{{\tilde
{\lambda}_{1}({\VPhi})}}

\newcommand{\lamtl}{{\tilde
{\lambda}_{2}({\VPhi})}}

\newcommand{\lmo}[1]{{\lambda_{1}({#1})}}
\newcommand{\lmt}[1]{{\lambda_{2}({#1})}}
\newcommand{\lmi}[1]{{\lambda_{i}({#1})}}

\newcommand{\lamo}{\lmo{\VPhi}}
\newcommand{\lamt}{\lmt{\VPhi}}
\newcommand{\lami}{\lmi{\VPhi}}

\newcommand{\Proofs}[2]{#1}

\begin{document}

\title{Pricing, competition  and market segmentation in ride hailing
}

\author{\IEEEauthorblockN{Tushar Shankar Walunj}
\IEEEauthorblockA{\textit{IEOR, IIT Bombay, India } \\
tusharwalunj@iitb.ac.in}
\and
\IEEEauthorblockN{Shiksha Singhal}
\IEEEauthorblockA{\textit{IEOR, IIT Bombay, India } \\
shiksha.singhal@iitb.ac.in}
\and
\IEEEauthorblockN{Veeraruna Kavitha}
\IEEEauthorblockA{\textit{IEOR, IIT Bombay, India } \\
vkavitha@iitb.ac.in}
\and
\IEEEauthorblockN{Jayakrishnan Nair}
\IEEEauthorblockA{\textit{EE, IIT Bombay, India } \\
jayakrishnan.nair@ee.iitb.ac.in}
\and
}

\maketitle

\begin{abstract}
We analyse a non-cooperative strategic game among two ride-hailing platforms, each of which is modeled as a two-sided queueing system, where drivers (with a certain patience level) are assumed to arrive according to a Poisson process at a fixed rate, while the arrival process of passengers is split across the two providers based on QoS considerations. We also consider two monopolistic scenarios: (i) each platform has half the market share, and (ii) the platforms merge into a single entity, serving the entire passenger base using their combined driver resources. The key novelty of our formulation is that the total market share is fixed across the platforms. The game thus captures the competition among the platforms over market share, which is modeled using two different Quality of Service (QoS) metrics: (i) probability of driver availability, and (ii) probability that an arriving passenger takes a ride. The objective of the platforms is to maximize the profit generated from matching drivers and passengers.

In each of the above settings, we analyse the equilibria associated with the game. Interestingly, under the second QoS metric, we show that for a certain range of parameters, no Nash equilibrium exists. Instead, we demonstrate a new solution concept called an \textit{equilibrium cycle}. Our results highlight the interplay between competition, cooperation, passenger-side price sensitivity, and passenger/driver arrival rates.

\end{abstract}


\section{Introduction}
There has been a significant interest in ride-hailing platforms like Uber, Lyft and Ola in recent years. These platforms have queues on two sides (passengers in need of a ride on one side and drivers on the other).
Such matching platforms have attained considerable scale. For example, the OLA company that provides a ride-hailing service in India has more than 1.5 million drivers across 250 cities. Passengers also often have the option of using a competing ride-hailing platform. Indeed, the friction associated with switching platforms is negligible from the standpoint of passengers. Thus, platforms must price their rides based on the interplay between price sensitivity on the passenger side, their own goal of revenue maximization, and also the crucial division of market share between competitors.

While there has been considerable interest in the performance evaluation of two-sided matching platforms in recent times (we provide a review of the related literature later in this section), most studies do not consider the impact of competition between platforms. This paper seeks to fill this gap. Specifically, we analyse a non-cooperative game between two ride-hailing platforms, where the platforms must compete for market share. The platforms compete via their pricing policies, the passenger base being both impatient as well as price sensitive. Additionally, for benchmarking, we also consider two monopolistic scenarios: (i) each platform has half the market share, and these market shares are insensitive to the pricing strategies of either platform, and (ii) the platforms merge into a single entity, serving the entire passenger base using their combined driver resources. 

Under the assumption of symmetric platforms, a single (geographical) zone of operation, and static (not state-dependent) pricing, we characterize the equilibria that emerge in the above scenarios, by approximating the payoff functions along a certain scaling regime (where driver impatience is diminishing). Comparing the different equilibria that emerge in these settings, our results highlight the impact of competition/cooperation between platforms, passenger price sensitivity, and passenger/driver arrival rates.

Our main contributions are as follows:

\begin{itemize}
    \item We model each platform as a BCMP network (\hspace{-.4mm}\cite{BCMP}), admitting a product form stationary distribution. This modeling approach is quite powerful, and allows for state-dependent pricing and driver impatience (it also permits  probabilistic routing across multiple zones, though this feature is not used in our analysis here). 
    \item We characterize the optimal pricing strategy of the platform under the monopolistic scenarios described above.
    \item We model market share bifurcation between the platforms in the form of a Wardrop equilibrium (\hspace{-.4mm}\cite{WE}) as in \cite{market_share}, where the passenger base splits in manner that seeks to equalize a certain QoS metric across the platforms. We consider two reasonable QoS metrics: the stationary probability of driver unavailability, and the stationary probability that an arriving passenger is not served (this latter metric also captures the scenario where the customer declines the ride  based on the high price quoted by the platform).
    \item Under each of the above passenger-side QoS metrics, we characterize the equilibria of the non-cooperative game between the providers. Interestingly, in the case of the second QoS metric, we show that under certain conditions, no Nash equilibrium exists. Instead, we demonstrate an \textit{equilibrium cycle}, where each platform has the incentive to set prices within a certain interval, though at each price (action) pair, at least one platform has the incentive to deviate to a different price within the same interval. 
    \item Finally, we compare the equilibria under the different settings  to highlight the impact of competition and cooperation between the platforms.
\end{itemize}


The remainder of this paper is organised as follows. After a brief survey of the related literature below, we describe our system model in Section~\ref{sec:model} and analyse the monopolistic setting in Section~\ref{sec:monopoly}. Competition under the driver availability QoS metric is analysed in Section~\ref{sec:pi0}, while the service availability metric (depends on the driver availability as well as the price quoted by the platform) is considered in Section~\ref{sec:overall_bp}. Finally, we report numerical experiments in Section~\ref{sec:comparison}, and conclude in Section~\ref{sec:conclusion}. Throughout the paper, references to the appendix point to the appendix of our technical report~\cite{TR}.

\subsection*{Literature Review}
Two sided queues have been considered in \cite{BCMP_routing},\cite{siva}-\cite{abandon},
where both customers and servers arrive over time, and wait
to be matched. 
Under a suitable fluid scaled limit
of this system, \cite{johari} showed that static pricing is sufficient to optimize the objective function, and dynamic pricing only
helps in improving the robustness of the system. The BCMP model in \cite{BCMP_routing} is almost similar to our model except the fact that we consider a open system (with varying number of drivers and passengers) with the possibility of drivers abandoning (and re-joining) the system. \cite{siva} considers a two-sided queueing system (with multiple zones and customer types) under joint pricing and matching controls. However, the possibility of the drivers re-joining the system has not been taken into consideration.

An extension of this two-sided queue model with multiple
types of servers and customers has also been considered
in \cite{uber}, \cite{kaplan}-\cite{ward}, where an additional bipartite matching
decision has to be made.  In \cite{kaplan}, limiting results of matching
rates between certain customer and server types with FCFS
scheduling have been analyzed. Without pricing, the optimal
matching to minimize the queueing cost has been analyzed
in \cite{ward} under a suitably scaled large system limit.  Most of these analyses are presented using a large
system scaling regime.

Queueing network models typically assume that the arrival process is exogenous. Contrasting system models have been considered in \cite{WE_app}-\cite{market_share} where the passenger arrival stream is split according to the well-known Wardrop Equilibrium, WE (~\cite{WE}). The authors in \cite{routing} consider that the  users choose the queue that will minimize their own expected delay. In \cite{market_share}, single-sided queueing system with multiple service providers has been considered; the objective of each of these providers is to maximize their market share. Further, the possibility of cooperation has been considered.

In contrast with the prior literature on two-sided queues, the important aspect of competition between platforms has not, to be best of our knowledge, been explored before.



\section{Model and Preliminaries}
\label{sec:model}

In this section, we describe our model for the interaction between two competing ride hailing platforms, and state some preliminary results. 

\subsection{Passenger arrivals}

We consider a system with a set $\mathcal{N} = \{1,2\}$ of independent ride-hailing platforms, each of which is modeled as a single-zone two-sided queuing system that matches drivers and passengers. Passengers arrive into the system as per a Poisson process with a rate $\Lambda$. This aggregate arrival process gets split between the two platforms based on QoS considerations (details in Subsection~\ref{sec:WE_QoS} below), such that
the passenger arrival process seen by platform $i$ is a Poisson process with rate~$\lambda_i$, where $\sum_{i=1}^2 \lambda_i = \Lambda$. 

When a passenger arrives into platform~$i,$ if there are no waiting drivers available, the passenger is immediately lost. On the other hand, if there are one or more waiting drivers, the passenger is quoted a price~$\phi_i(\niw) \in [0,\phi_h],$ where~$n_i$ denotes the number of waiting drivers on platform~$i,$ and $\phi_h$ denotes the maximum price the platform can charge. The customer accepts this price (and immediately begins her ride) with probability~$f(\phi_i(\niw));$ with probability $1-f(\phi_i(\niw)),$ the passenger rejects the offer and leaves the system (without taking the ride). Note that the function~$f$ captures the price sensitivity of the passenger base. We make the following assumptions on the function~$f.$
\begin{enumerate}[{\bf A}.1]
\item  $f(\cdot)$ is a strictly concave, strictly decreasing and differentiable function.
\item $0 < f(\phi) \le 1$ for all $\phi \in [0,\phi_h]$ and $f(0) =1$.
\end{enumerate}

\subsection{Driver behavior}

Each platform has a pool of dedicated drivers. Dedicated drivers of platform~$i \in \mathcal{N}$ arrive into the system according to a Poisson process of rate $\eta_i$. The drivers wait in an FCFS queue to serve arriving passengers. Recall that the number of waiting drivers in this queue on platform~$i$ is denoted by $\niw$. 

When an arriving passenger is matched with the head of line driver, a ride commences, and is assumed to have a duration that is exponentially distributed with rate $\nu_i.$ At the end of this ride, the driver rejoins the queue of waiting drivers (and therefore becomes available for another ride) with probability~$p_i;$ with probability~$1-p_i,$ the driver leaves the system.

\begin{figure}[t]
     \centering \includegraphics[scale=0.25]{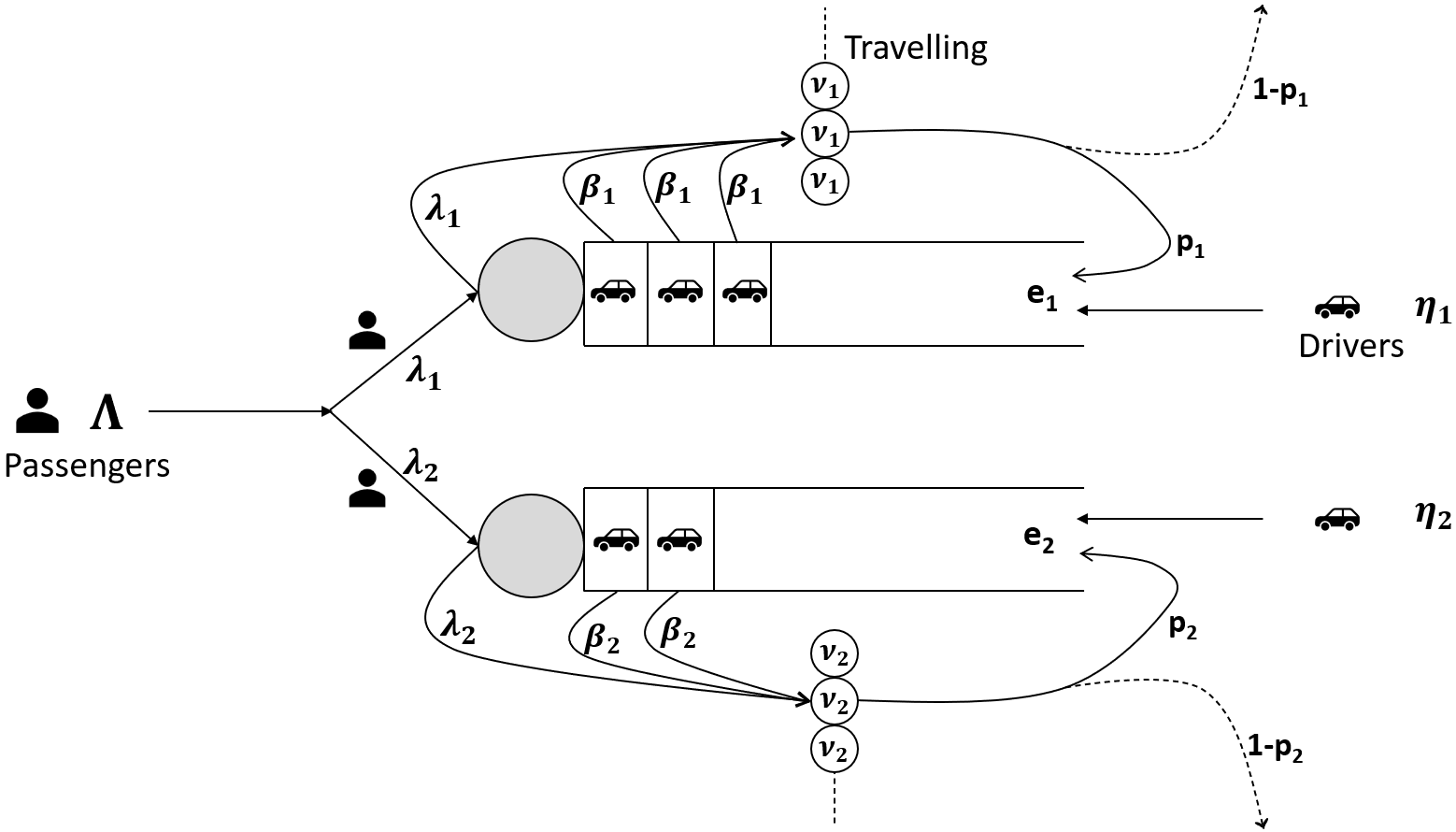}
       \caption{Depiction of system model}
       \label{fig:model}
 \end{figure}

Additionally, we also model driver abandonment from the waiting queue. Specifically, each waiting driver independently abandons after an exponentially distributed duration of rate~$\beta_i.$ This abandonment might capture a local (off platform) ride taken up by the driver, or simply a break triggered by impatience. The duration of this ride/break post abandonment is also assumed to be exponentially distributed with rate~$\nu_i$ (when drivers go for local rides it is reasonable to assume that the local ride duration has the same distribution, as the local ride is also in the same area); at the end of this duration, the driver rejoins the queue with probability~$p_i$ and leaves the system altogether with probability~$1-p_i.$ Figure~\ref{fig:model} presents a pictorial depiction of our system model. 

\subsection{BCMP modeling of each platform}

Under the aforementioned model, each platform can be described via a continuous time (Markovian) BCMP network~\cite{BCMP}. Formally, the state of platform $i$ at time $t$ is given by the tuple $Z_{i,t} = (N_{i,t},R_{i,t})$, where $N_{i,t}$ is the number of waiting drivers, and $R_{i,t}$ is the number of drivers that are in the system but unavailable to be matched with passengers (because of being in a ride, or on a break). Realizations of the state of platform~$i$ are represented as $(\niw, \nir)$, and state space corresponding to platform~$i$ is given by $S_i = \{(\niw, \nir)\colon \niw, \nir \in \mathbb{Z}_+\}$. 

The BCMP model for each platform consists of two `service stations': (i) service station 1 (SS1) is the queue of waiting drivers (the occupancy of this queue is the first dimension of the state), and (ii) 
service station 2 (SS2) is the `queue' of drivers on a ride/break (the occupancy of this queue is the second dimension of the state). SS1 is modeled as as a single server queue having a state dependent service rate. Specifically, the departure rate from SS1 is $\lambda_i + \niw \beta_i.$ On the other hand, SS2 is modeled as an infinite server queue, having exponential service times of rate~$\nu_i.$ SS1 sees exogenous Poisson arrivals at rate~$\eta_i,$ and departures from SS1 become arrivals into SS2. Finally, departures from SS2 exit the system with probability~$1-p_i,$ and join SS1 with probability~$p_i;$ see Figure~\ref{fig:model}.

Next, we describe the steady state distribution corresponding to the above BCMP network (associated with platform~$i$). Define $e_i := \frac{\eta_i}{1-p_i},$ which is easily seen to be the effective arrival rate seen by the service stations SS1 and SS2. The following lemma follows from \cite{BCMP} (see Sections~3.2 and~5).  
\begin{lemma} \label{lem_BCMP_Product_form}
The steady state probability of state $s = (\niw, \nir)\in S_i$ is given by:
\begin{align*}
\pi_i(s) &= C_i
	\left[\frac{{e_i}^{\niw}}{\prod_{a=1}^{\niw} \left(\lambda_i f(\phi_i(a))+ a \beta_i\right)}\right]
	\left[\left({\frac{e_i}{\nu_i}}\right)^{\nir} \left(\frac{1}{\nir!}\right) \right], \\
C_i^{-1} 
&= \sum_{n=0}^{\infty} \left[\frac{e_i^{n}}{\prod_{a=1}^{n}(\lambda_i f(\phi_i(a))+ a \beta_i)}\right]
\exp\left({\frac{e_i}{\nu_i}}\right).
\end{align*}
Here, $C_i$ is the normalizing constant (with the convention that $\prod_{a}^{k} (\cdot) = 1$ when $ a> k$). \eop
\end{lemma}

\subsection{Passenger split across platforms: Wardrop Equilibrium}
\label{sec:WE_QoS}

We model the split of the aggregate passenger arrival rate~$\Lambda$ into the system across the two platforms (recall that platform~$i$ sees passenger arrivals as per a Poisson process with rate~$\lambda_i$) as a Wardrop equilibrium (WE) based on Quality of Service (QoS). In particular, we consider two different QoS metrics in our analysis. However, before describing these formally, we first define the Wardrop split of the passenger arrival rate in terms of a generic QoS metric~$Q.$

Let $Q_i(x)$ denote the QoS of platform $i$ when the passengers arrive at rate~$x.$ Note that~$Q_i$ will, in general, also depend on the pricing policy~$\phi_i := (\phi_i(n),\ n \in \mathbb{N})$ employed by platform~$i,$ though this dependence is suppressed for simplicity. The Wardrop split $(\lambda_1,\lambda_2)$ under the price policy $\VPhi = (\phi_1, \phi_2)$ is then defined as:
\begin{equation}\label{eqn_existence_uniq_WE_main}
 \lamo \in \arg \min_{\lambda \in [0,\Lambda] } \left(Q_1  (\lambda) - Q_2 (\Lambda-\lambda) \right )^2, \ \lamt = \Lambda - \lamo.  
\end{equation}
We address the uniqueness of the Wardrop split in Lemma~\ref{lem_exist_unique_WE} below. Having defined the Wardrop split in generic terms, we now define the passenger QoS metrics we consider. Note that passengers can leave the system without taking a ride either because of driver unavailability, or because the price quoted was too high.


Let $\DA_i$ be the stationary probability of zero waiting drivers on platform~$i$. From Lemma \ref{lem_BCMP_Product_form}, 
\vspace{-2mm}
 \begin{equation}
\label{eqn_pi0_expression}
\DA_i = \sum_{\nir = 0}^{\infty} \pi_i((0, \nir)) \ = \left(\sum_{\niw=0}^{\infty} \left[\frac{(e_i)^{\niw}}{\prod_{a=1}^{\niw}(\lambda_i f(\phi_i) + a \beta_i)}\right]\right)^{-1} \hspace{-3mm}.
\end{equation}
By PASTA (\cite{pasta}), it follows that $\DA_i$ is also the long run fraction of passengers who find no available driver upon arrival into platform~$i$ (and thus leave the system).

Let $\PB_i$ be the long run fraction of passengers who leave platform~$i$ without taking a ride (due to driver unavailability or a high price). Using PASTA (\cite{pasta}) again, 
\begin{align}
\PB_i &= \DA_i + \sum_{\niw=1}^{\infty} \sum_{\nir=0}^{\infty} \left[ (1 - f(\phi_i(\niw))) \pi_i((\niw, \nir))\right]. \label{eqn_expression_of_PBi}
\end{align}

We use the $\DA_i $ and $\PB_i$ as our QoS metrics; note that these are both functions of $\lambda_i$. Our analysis for the case~$Q_i = \DA_i$ is presented in Section~\ref{sec:pi0}, while the case~$Q_i = \PB_i$ is addressed in Section~\ref{sec:overall_bp}. Next, we establish existence and uniqueness of the Wardrop equilibrium under these metrics (proof in Appendix~A); we first prove existence and uniqueness under a certain assumption (Assumption~{\bf A}.0), and them show that our metrics of interest satisfy the assumption.

\ignore{
We use these two metrics for QoS and study the above mentioned model. Let $Q_i(\lambda_i)$ denote any given quality of service  which influences the Wardrop equilibrium.
\begin{enumerate}
    \item \textbf{$(Q_i(\lambda_i) = \DA_i$:} Passenger blocking due to driver unavailability,
    \item \textbf{$(Q_i(\lambda_i) = \PB_i)$:} Passenger blocking due to driver unavailability and unacceptable price.
\end{enumerate}
We have the following result under the following assumption
 \begin{enumerate}[{\bf A}.0]
     \item The QoS function $Q_i$ for each $i$ is continuous and strictly monotone in $\lambda_i$; both the functions are either increasing or decreasing.
 \end{enumerate}
} 

 \begin{enumerate}[{\bf A}.0]
     \item The QoS function $Q_i$ for each $i$ is continuous and strictly monotone in $\lambda_i$; both the functions are either increasing or decreasing.
 \end{enumerate}
 
\begin{lemma}
\label{lem_exist_unique_WE}
Consider any QoS metric satisfying Assumption~{\bf A}.0. Given any price policy $\VPhi = (\phi_1, \phi_2)$, there exists a unique Wardrop Equilibrium $(\lamo,\lamt).$

Finally, each of the QoS metrics $\DA_i$ and $\PB_i$ satisfies {\bf A}.0, so long as $\beta_i > 0.$ \eop
\end{lemma}

Given the existence and uniqueness of WE, we next define the platform utility functions. This defines a non-cooperative game between the two platforms.

\subsection{Platform Utilities}
\label{ssec:strategic_game}

We treat the action of each platform~$i$ to be its pricing policy~$\phi_i = (\phi_i(n),\ n \in \mathbb{N}).$ We define the utility of platform~$i$ as the (almost sure) rate at which it derives revenue from matching drivers with passengers, denoted by $\MR_i.$ Note that $\MR_i$ depends on the price (action) profile~$\VPhi = (\phi_1, \phi_2)$ and other parameters; the pricing policy of each platform influences the other's utility via the Wardrop split that determines the market shares of both platforms.
\ignore{
Consider platforms as described above. The system consists of two competing platforms, competing for the maximum profit (the matching revenue in our case).
Consider a strategic form game,
\vspace{-.5mm}
\begin{equation}\label{eqn_game}
\left< \mathcal{N} = \{1,2\}, \left(\phi_i \in [0, \phi_h]\right)_{i \in \mathcal{N}}, \left(\MR_i \in [0, \infty) \right)_{i \in \mathcal{N}} \right>
\end{equation}
where the player $i$ has a utility as the matching revenue $\MR_i$ and $\phi_h$ is the maximum price that any platform can offer; observe $\MR_i$ depends upon the Wardrop split $(\lamo, \lamt)$ for any price vector $\VPhi$.
} 
In the following, we derive the matching revenue rate (often referred to simply as the \textit{matching revenue}, or MR) of each platform in terms of the Wardrop split $(\lambda_1,\lambda_2)$ (proof in Appendix A).
\begin{lemma} \label{lem_mr_derivation}
The matching revenue rate of platform $i$ 
is given by
\begin{equation} \label{eqn_MR_expression}
\MR_i =  \sum_{s_i: \niw \neq 0} \lambda_i f(\phi_i(\niw)) \phi_i(\niw)\pi_i(s_i).
\end{equation} 
\end{lemma}
\noindent

\subsection{Static prices and other simplifications}
The goal of this paper is to analyse the game defined above between the two platforms. Since the analysis under dynamic pricing policies (where a platform's price varies with the number of waiting drivers) is quite complex, we restrict attention to static pricing policies in this paper. Indeed, even the analysis under static pricing is non-trivial and highlights several novel insights. Note that under static pricing, the action space for each platform is $[0,\phi_h].$



Further, for simplicity of presentation, we consider symmetric platforms, i.e., we assume $e_1 = e_2 = e$ and  $\beta_1=\beta_2 = \beta.$\footnote{Our results can be generalized to the case where the $e_i$ are distinct.}

Finally, we approximate the utility functions of the platforms by letting~$\beta \downarrow 0.$ While this limiting regime enables an explicit characterization of the equilibria of the non-cooperative game under consideration, it is also well motivated. Indeed, $\beta \approx 0$ means that driver abandonment times are stochastically much larger than ride durations and passenger interarrival times. This is reasonable in many practical scenarios, particularly when drivers are `tied' to their respective platforms, and are available for rides for long durations (say an 8-10 hour work shift) at a stretch.

We conclude this section with a result characterizing the limit of the MR of any platform as $\beta \to 0$ (proof in Appendix~A).\footnote{Throughout, in our notation, we emphasize functional dependence on parameters of interest only as and when required.}


\begin{lemma} \label{lem_approx_MR}
Consider static policy profile $(\phi_1,\phi_2)$. Suppose that the passenger arrival rates $\lambda_i (\beta) \to {\lambda}_i$ as $\beta \to 0$. Then the matching revenue $\MR_i = \MR_i (\phi_i;\lambda_i(\beta), \beta),$ given by Lemma~\ref{lem_mr_derivation}, converges to $\MRl_i$ as $\beta \to 0,$ where
\begin{align}
\label{eq_approx_MR}
\MRl_i(\phi_i) &=  \begin{cases}
        \eiw \phi_i & \text{ if }  \left(\frac{\eiw}{\lambda_i f(\phi_i)}\right) < 1, \\
        \lambda_i f(\phi_i) \phi_i & \text{ else.}
    \end{cases} 
\end{align}
Along similar lines, we have  $\DA_i(\phi; \lambda_i(\beta), \beta)  \to \DAl_i (\phi; \lambda_i)$, where
\begin{align}
\label{eq_pi0_approx_lem}
\DAl_i(\phi_i; \lambda_i) &=
        \left(1-\frac{\eiw}{\lambda_i f(\phi_i)}\right)^+.
\end{align}
\end{lemma}

By virtue of the above result, we derive approximate equilibria of the actual system when $\beta$ is small, by analysing those of a `limit system,' which is obtained by letting~$\beta \to 0.$ Throughout, we refer the system (i.e., functions) obtained from the limit $\beta \to 0$ as the limit system. In fact, we also provide a formal justification of our conclusions from this limit system by establishing convergence of relevant optimizers and zeros, as and when required. 


\ignore{\begin{table}[t]
    \centering
    \begin{tabular}{|c|c|}
   \hline 
  $  m (\phi) = \Lambda f(\phi) \phi - e \phi $  Matching Revenue - for $\bar{\phi}> \phi > \underline{\phi}$ &   $  \phi^*_m $  \\ \hline
         $ {M} (\phi) =  \Lambda f(\phi) \phi$  
      Matching Revenue - for $\phi > \bar{\phi}$    &  $ \phi^*_M$   \\ \hline 
     $\phiunderbar$  &  $f^{-1} (2e/\Lambda)$  \\ \hline  
      $\phibar$ &  $f^{-1} (e/\Lambda)$  \\ \hline
   $ d(\phi) =   2f(\phi) + \phi f'(\phi)$ & $\phi^o_d$ \\ \hline
    \end{tabular}
    \caption{Caption}
    \label{tab:notations}
\end{table}}

\section{Monopoly}
\label{sec:monopoly}

In this section, we consider the scenario where a single platform operates alone (with no competitors), with an incoming passenger arrival rate half that of original system, i.e.,  $\nicefrac{\Lambda}{2}$. The remaining details (like passenger price sensitivity, driver behavior, etc.) are as described in Section~\ref{sec:model}.
This is equivalent to the special case of our model where the passengers arrival rate is split equally between the two platforms, independent of the platforms' pricing policies.\footnote{This special case arises by taking~$Q_i(\lambda_i) \equiv 1.$} In this case, the non-cooperative game reduces to separate (and identical) optimization problems for the platforms to maximize their revenue rate; this optimization is analysed in this section.

The platform offers price $\phi$ to the incoming passenger, if there are waiting drivers  and then  the passenger accepts    the price   with probability $f(\phi)$. The matching revenue derived by the platform is given by Lemma~\ref{lem_mr_derivation} with $\lambda = \Lambda/2$ and the static policy $\phi(n)=\phi$ for all $n$ (see \eqref{eqn_MR_expression}): 
\begin{equation}
\label{eqn_mr_monopoly}
\MR(\phi; \beta) = \frac{\Lambda}{2} f(\phi) \phi (1 - \DA(\phi; \beta)) 
\end{equation}
This can be approximated using Lemma~\ref{lem_approx_MR} when $\beta$ is close to zero. Formally, we consider the approximate platform utility corresponding to the limit system:
\begin{align*}
    \MRl(\phi) 
    = \begin{cases}
        e \phi & \text{ if }  \left(\frac{e}{\frac{\Lambda}{2} f(\phi)}\right) < 1 \\
        \frac{\Lambda}{2} f(\phi) \phi & \text{ else}
    \end{cases}
\end{align*} 

The optimal pricing strategy for the monopolistic platform, seeking to maximize $\MRl(\phi),$ is characterized as follows (proof in Appendix~A).

\begin{thm}\label{thm_monopoly_opt}
The optimizer of the MR in the limit system $\MRl(\phi)$ \eqref{eq_approx_MR} is given by,
\begin{align}
\label{Eqn_Mphi_PhistarMPhi}
    \phi^* &= \max\{\phifp, \phiunderbar\}, \mbox{ where } 
    \phiunderbar := f^{-1}\left(\frac{2e}{\Lambda}\right), \mbox{ and } \\ \nonumber 
    \phifp &:= \arg \max_{\phi} \P (\phi) \mbox{, with, }  \P (\phi) =  f(\phi) \phi.  \hspace{18mm}
\end{align} 
Moreover, for any sequence of optimal prices corresponding  to a  given sequence $\beta_n \to 0$, there exists a sub-sequence that converges to the unique optimal price of the limit system. \eop
\end{thm}

Theorem~\ref{thm_monopoly_opt} characterizes the optimal policy of the monopolistic platform in the limit system, and also justifies this approximation when~$\beta$ is small. In following sections, we will find it instructive to compare the equilibrium pricing policy of each platform to the monopolistic optimal policy \eqref{Eqn_Mphi_PhistarMPhi}, to shed light on the impact of passenger-side churn on platforms' pricing policies.

\section{Duopoly driven by  Driver Availability}
\label{sec:pi0}

We now return to the platform duopoly model introduced in Section~\ref{sec:model}.
Specifically, in this section, we assume that the passengers are primarily sensitive to the non-availability of drivers, i.e., to the probabilities $\{\DA_i\}$. Accordingly, the QoS metric $Q_i$ is taken to be $\DA_i,$ in order to define the Wardrop split~\eqref{eqn_existence_uniq_WE_main} of passenger arrival rate across the two platforms.

This results in a non-cooperative 
game between the two platforms. We now seek to characterize the Nash equilibria of this game in the limit system. 
%
%
We begin with a characterization of the WE (proof in Appendix~A). 
\begin{lemma} \label{lem_pio_WE}
For a given price vector $\VPhi, \beta>0$ and passenger response function $f$, the WE under the QoS metric $\DA_i$ is given by
\begin{equation*}
  (\lamo, \lamt)   = 
 \left(\frac{\Lambda f(\phi_2)}{f(\phi_1) + f(\phi_2)}, \frac{\Lambda f(\phi_1)}{f(\phi_1) + f(\phi_2)} \right).   
\end{equation*}
\end{lemma}

The above lemma provides the unique WE for any given price-vector/strategy profile $\VPhi,$ and for any~$\beta > 0$ (i.e., we have not needed to appeal to the limit system). Indeed, note that the Wardrop split is insensitive to the value of~$\beta$ under QoS metric~$\DA_i.$ Now, the utility of player $i$, i.e., its matching revenue $\MR_i (\VPhi) = \MR_i(\phi_i; \lami)$, can be approximated by  $\MRl_i(\VPhi)$ in the limit system; this approximation is obtained from~\eqref{eq_approx_MR} after replacing $\lambda_i $ of Lemma~\ref{lem_approx_MR} by $\lami$ of Lemma~\ref{lem_pio_WE}. We now derive the NE of this `limit game' (proof in Appendix~A).

\ignore{
\begin{thm} \label{thm_pio_NES}
If there exists a $ \phido $ such that $d(\phido)= 0$, and
\begin{enumerate}[(a)]
\item  $\phiunderbar \in [0, \phi_h]$, then $(\phi^*, \phi^*)$ where $\phi^* = \max \{ \phiunderbar,\phido\}$ is a Nash equilibrium (NE),
    \item $ \phiunderbar \notin [0, \phi_h]$ {\color{red}i.e., $e > \Lambda$}, then $(\phido, \phido)$ is a NE,
\end{enumerate} 
else 
$(\phi_h, \phi_h)$ is a NE. \eop
\end{thm}}

\begin{thm} 
\label{thm_pio_NES}
 Define $d(\phi) := (2f(\phi) + \phi f'(\phi))$. Then\\
i) $(\phiunderbar, \phiunderbar)$ is a Nash equilibrium (NE) if and only if  $d (\phiunderbar) \le  0$;\\
ii) If $d (\phiunderbar) > 0$ and 
$d(\phi_h) \le 0,$ 
then  $ (\phido, \phido)$ is a NE, where $d(\phido) = 0$; \\ 
iii) if not the first two cases, then $(\phi_h,\phi_h)$ is a NE. \eop
\end{thm}
Theorem~\ref{thm_pio_NES} can be interpreted as follows. Given our assumptions on the price sensitivity function~$f,$ it is easy to see that~$d$ is strictly decreasing. Case~$(iii)$ above arises when~$d(\phi_h) >0;$ this means passengers are (relatively) price insensitive. In this case, the NE corresponds to both platforms charging the maximum permissible price. On the other hand, Cases~$(i)$ and~$(ii)$ arise when $d(\phi_h) \leq 0,$ i.e., when passengers are (relatively) price sensitive. In particular, the former case arises when $e/\Lambda$ is small (i.e., the passenger arrival rate is large relative to the rate at which drivers are available), while the latter case arises when $e/\Lambda$ is large (i.e., the passenger arrival rate is small relative to the rate at which drivers are available). The equilibrium price decreases as a function of $e/\Lambda;$ in other words, when the rate at which passengers enter the system decreases, the equilibrium price quoted by the platforms also decreases. We compare the values of the equilibrium prices here to the monopolistic optimal price (characterized in Section~\ref{sec:model}) in Section~\ref{sec:comparison}. 

Finally, a (partial) justification of the equilibria under the limit system can be provided as follows. The  best response (BR) of   player $i$ against any given $\phi_{-i}$ of the other player, converges to the corresponding BR in limit system in the same sense (along a sub-sequence) as in Theorem~\ref{thm_monopoly_opt}; the proof follows using exactly the same logic as that in Theorem~\ref{thm_monopoly_opt}.

\section{Duopoly  driven by Price \& Driver Availability}
\label{sec:overall_bp}

Contrary to the previous section, we now consider the scenario where passengers are sensitive to overall likelihood of finding a ride on arrival; note that the ride may not materialise either due to driver unavailability, or due to the price being prohibitive. Formally, we choose the QoS metric to be the combined blocking probability $\{\PB_i\}.$ 
Thus Wardrop split \eqref{eqn_existence_uniq_WE_main}  ensures the difference in the combined blocking probabilities $\{\PB_i\}$ of the two platforms is minimized.
 
 The combined blocking probability of platform $i$ when it uses static price policy $\phi_i$ and when passengers arrive at rate $\lambda_i$ is given by (see \eqref{eqn_expression_of_PBi}),
 
\vspace{-4mm}

{\small \begin{align}
\label{Eq_overall_exact}
    \PB_i(\phi_i) &= 
     \DA_i  + (1-f(\phi_i))(1-\DA_i)   =  \DA_i f(\phi_i) +(1-f(\phi_i)) ,
    \end{align}}
    \noindent which can be approximated by considering the limit system as follows (using Lemma \ref{lem_approx_MR}),
    \begin{equation}
    \label{Eq_overall_aproxx}
   \PBl_i(\phi_i; \lambda_i)  = \begin{cases}
            1 - \frac{\eiw}{\lambda_i} & \text{ if }  \left(\frac{\eiw}{\lambda_i f(\phi_i)}\right) < 1, \\
            1 - f(\phi_i) & \text{ else.}
        \end{cases}
\end{equation}



As before, the matching revenue of platform $i$ after WE split, when the two platforms operate with static price policies $\VPhi$, equals $\MR_i (\VPhi) = \MR_i (\VPhi; \lami)$, where  $\lami$  is now defined using QoS $\PB$.
%
As in the previous sections, we derive approximate platform utilities when $\beta \to 0$. However, we require a second level of approximation in this case; we first approximate the WE and then the corresponding MR to obtain  the utility functions in the limit system.


\begin{thm}
\label{thm_WE_MR_PB}
Fix $\VPhi, i$.
Let the (unique) WE be represented by $W(\beta) =  \lambda_{1}(\VPhi; \beta) $ for any $\beta>0$. For $\beta=0$, define $W(0) = \lamol$ as defined in Table \ref{tab:WE_MR_PB}. Then, the mapping 
$\beta \mapsto W(\beta) $ is continuous on the interval $[0,\infty)$. Further the same is the case with matching revenue function
$\beta \mapsto \MR_i(\VPhi; \beta )$, when we set $  \MR_i (\VPhi;0) := \MRl_i (\VPhi)$, where $\MRl_i$ is defined in Table \ref{tab:WE_MR_PB}. \eop

\end{thm}

\begin{table}[h]
\begin{center}
\begin{tabular}{|c|c|c|c|c|}
\hline
&&&& \vspace{-2mm} \\  
& $\phi_1 < \phiunderbar$ & $\phi_1 \in [\phiunderbar, \phibar)$                        & $\phi_1 \in [\phibar, \phi_h]$   & $\phi_1 \in [\phiunderbar, \phi_h]$ \\  
&  $\phi_1 \ge  \phi_2$    &  $\phi_1 > \phi_2$   & $\phi_1 > \phi_2$ & $\phi_1 = \phi_2$  \\ \hline 

$\lamol$                                                    & $\frac{\Lambda}{2}$            & {\scriptsize$\left(\Lambda - \frac{\eow}{f(\phi_1)}\right)$}           & $0$    & { $\frac{\Lambda}{2}  $}                                 \\ \hline
$\MRl_1(\VPhi)$                                                    & $\eow\phi_1$                   & $m(\phi_1)$                    & $0$   & {\scriptsize$\frac{\Lambda}{2}f(\phi_1)\phi_1$}                                     \\ \hline
$\MRl_2(\VPhi)$                                                    & $\etw \phi_2$                  & $\etw \phi_2$                                               & {\scriptsize  $\min\{\Lambda f(\phi_2), \etw \} \phi_2$}  &  {\scriptsize$\frac{\Lambda}{2}f(\phi_2)\phi_2$} \\ \hline
\end{tabular}

\end{center}
\caption{Limit system under QoS $\PB$, $\phiunderbar = f^{-1}\left(\frac{2e}{\Lambda} \right)$, $\phibar = f^{-1}\left(\frac{e}{\Lambda} \right)$ (set $\phiunderbar,\phibar = 0 $ when not in range) and $m(\phi) = (\Lambda f(\phi)-\eow)\phi$ \label{tab:WE_MR_PB}
}
\end{table}

By virtue of the above theorem, the MR of the actual system can be approximated by that of the limit system as defined in Table \ref{tab:WE_MR_PB}, when $\beta$ is close to zero.  Once again, by this continuity, the  best response (BR) of   player $i$ against any given $\phi_{-i}$ of the other player, converges to the corresponding BR in the limit system (along a sub-sequence) as in Theorem \ref{thm_monopoly_opt}.  

So for similar reasons as before,  we henceforth consider the limit system to  obtain the corresponding equilibrium. 
However, there are some interesting contrasts between the results in this section and those in Section~\ref{sec:pi0}. The most important of them being that for a certain range of parameters, we do not have a Nash Equilibrium (NE). This is possibly because of the discontinuous utilities that players in such a duopoly derive. Basically given  any  price vector, one can find an alternate price  that yields them strictly better revenue upon unilateral deviation;  interestingly, in our case the choices of better prices seem to be confined to an interval. 
This motivates us to propose \textit{an alternate equilibrium concept, which we term as an equilibrium cycle. We show that an equilibrium cycle indeed exists for such a duopoly.} We begin with deriving the conditions under which classical Nash Equilibrium (NE) exists.

\begin{thm}
\label{thm_PB_NE_phiunderbar}
Assume $e < \Lambda.$ Then $(\phiunderbar, \phiunderbar)$ is a NE
if $ \phimr \leq \phiunderbar \leq \phi_h$, where $\phimr$ is the unique maximizer of $m(\phi) = (\Lambda f(\phi)- e) \phi$ over $[0,\phi_h].$ \eop
\end{thm}

In the regime of parameters that do not satisfy the hypothesis of the above theorem, we don't have classical NE. We now state the definition of proposed equilibrium cycle and show the existence of the same. 

\begin{defn} \label{def_equi_cycle}
{\bf Equilibrium Cycle:} A closed interval $[a, b]$ is called an equilibrium cycle, if 

(i) for any $ i \in \mathcal{N}$, $\phi_{-i} \in [a,b]$ and $ \phi_i \in [0,\phi_h] \setminus [a, b] $, there exists a $\phi_i' \in [a,b]$, we have $\MR_i(\phi_i', \phi_{-i}) >  \MR_i(\VPhi)$, and

(ii) for any price vector $\VPhi \in [a,b]^2$, there exists $ i \in \mathcal{N}$ and $ \phi_i' \in [a,b]$ such that $\MR_i(\phi_i', \phi_{-i}) > \MR_i(\VPhi)$.
\end{defn}
The first condition above establishes the `stability' of the interval~$[a,b];$ if any player has an action in this interval, the other player is also incentivized to play an action in the same interval. The second condition establishes the `cyclicity' of the same interval; if both players play any actions within the interval~$[a,b]$, at least one player has an incentive to deviate to a different action within the same interval. We show below that we indeed have an equilibrium cycle $[\eL, \eU]$ as given by the following (recall that $\phimr$ is  defined in Theorem \ref{thm_PB_NE_phiunderbar}):
$$
\left[\eL, \eU \right]  =  \left[ \left(\frac{\Lambda f(\phimr)}{e} - 1\right)\phimr,\phimr \right] = 
\left[ \ \frac{m(\eU)}{e}  ,\phimr \right]  
.
$$

\begin{thm}
\label{Thm_PB_EC} Assume $e < \Lambda.$
If $\phiunderbar < \phimr \le \min\{\phi_h, \phibar\}$ then $\left[\eL, \eU \right]$ is an equilibrium cycle. \eop
\end{thm}

The proof of the above theorem is in Appendix A.
For the remaining range of parameters, when the duopoly is driven by price as well as driver availability, there is neither a NE nor an equilibrium cycle. The market share is small and the competition is high, so none of the platforms find it beneficial to quote high prices. It can also be shown that $(0,0)$ is also not a NE. However, we do have an $\epsilon$-NE suggesting that players choose a price really close to zero in such highly competitive scenarios.

\begin{thm}\label{thm_epsilon_NE}
Consider the case with $e \ge \Lambda$. For any $\epsilon > 0$, if possible, choose $0<\delta \le \phi_h$ such that   $\sup_{\phi \leq \delta} \Lambda f(\phi)\phi < \epsilon$. Then $(\delta, \delta)$ is an $\epsilon$-Nash equilibrium. \eop
\end{thm}

Jointly, Theorems~\ref{thm_PB_NE_phiunderbar}--\ref{thm_epsilon_NE} characterize the equilibrium behavior in the limit system. We compare these equilibria with those derived in Section~\ref{sec:pi0} as well as the monopolistic optimal pricing in Section~\ref{sec:model} in Section~\ref{sec:comparison}.

\input{coop.tex}

\begin{figure}[t]
     \centering
     \begin{minipage}{3.7cm}
     \includegraphics[trim = {0cm 6cm 1cm 6cm}, clip, width = 3.8cm, height = 4cm]{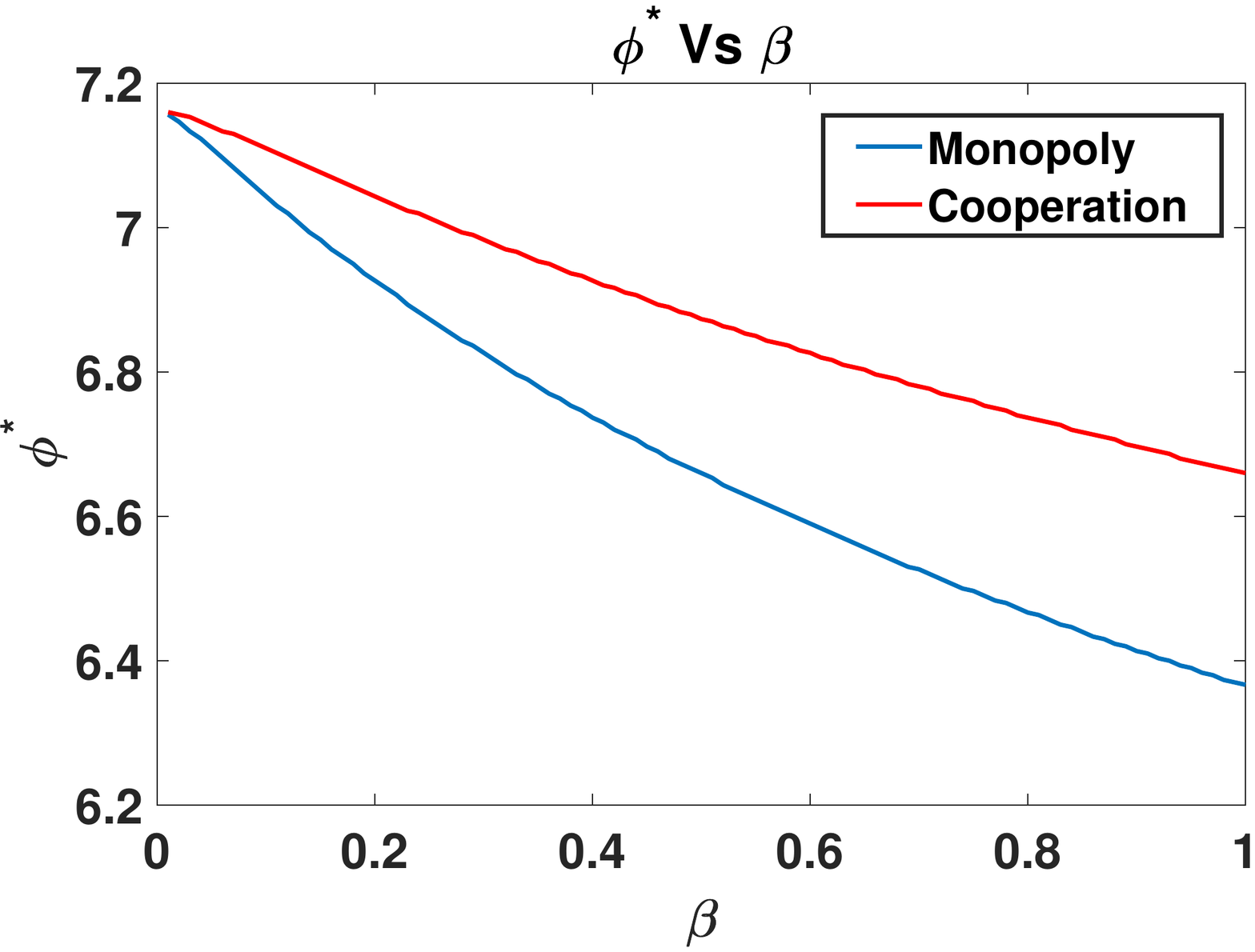}
             \end{minipage}
        \hspace{3mm}     \begin{minipage}{3.7cm}\includegraphics[trim = {0cm 6cm 1cm 6cm}, clip, width = 4cm, height = 4cm]{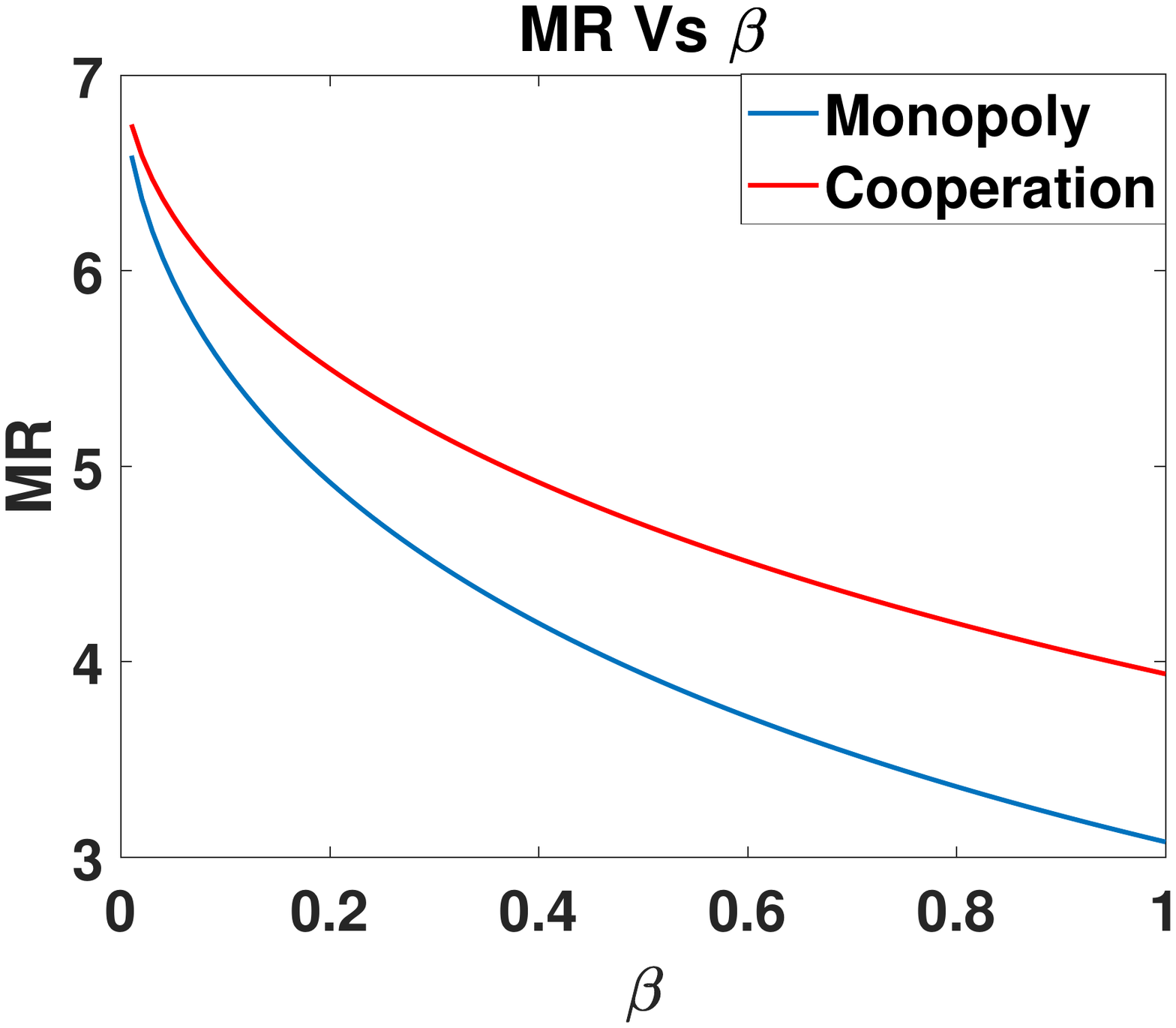}
             \end{minipage}
       \caption{Optimizers and optimal MR versus $\beta$; theoretical $\phi^*$ for the limit system is 7.1436; $e = 1$, $\Lambda = 7$, $f(\phi) = 1- a \phi$, with $a = 0.1$}
       \label{fig:monopoly_coop}
 \end{figure}

\begin{figure}[h]
    \centering
    \hspace{-15mm}
    \begin{minipage}{3.7cm}
    \includegraphics[width = 4.4cm, height = 4cm]{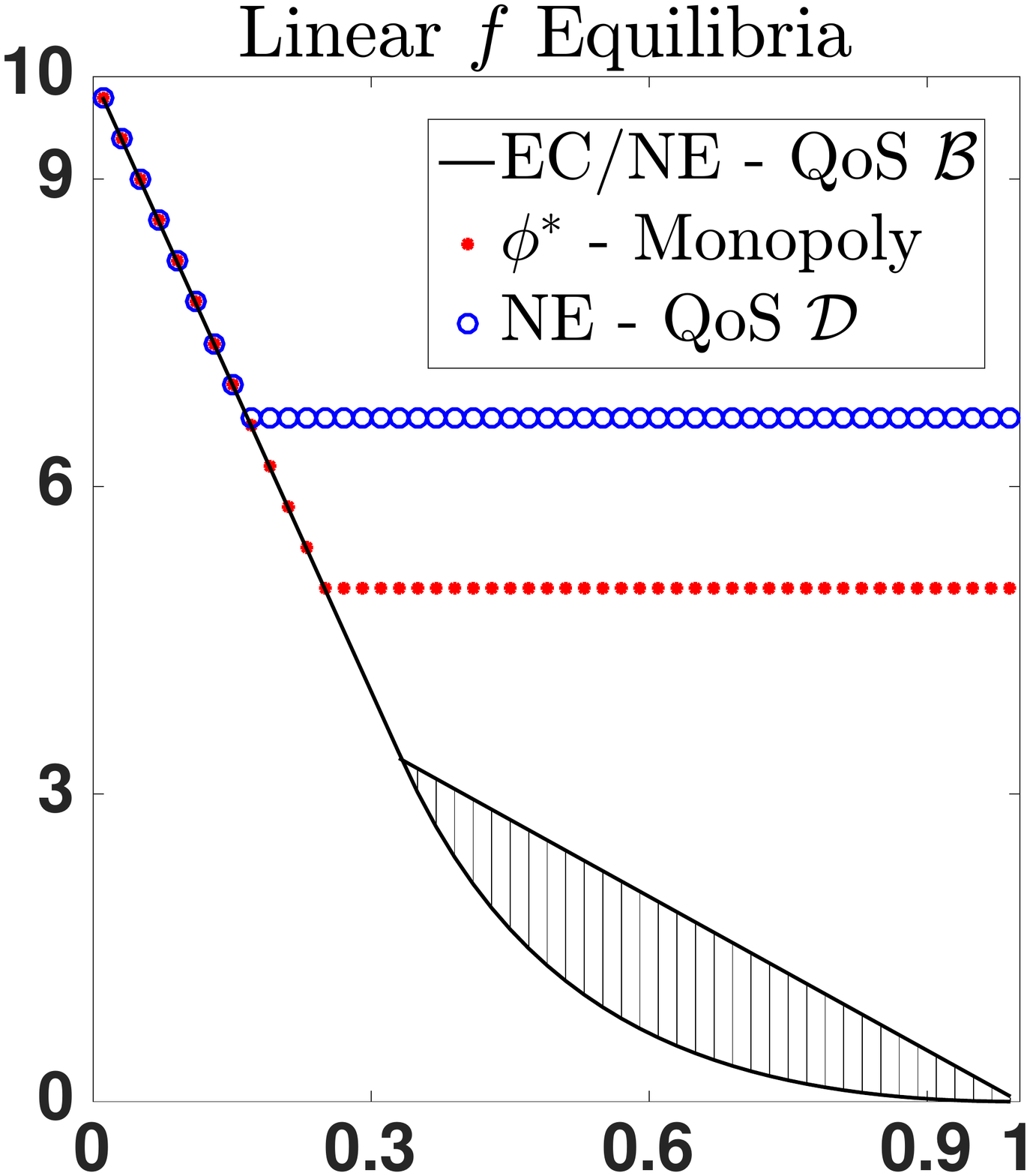}

        \end{minipage}
        \hspace{6mm}
      \begin{minipage}{3.7cm}
      \vspace{-2mm}
         \includegraphics[width = 4.5cm, height = 4cm]{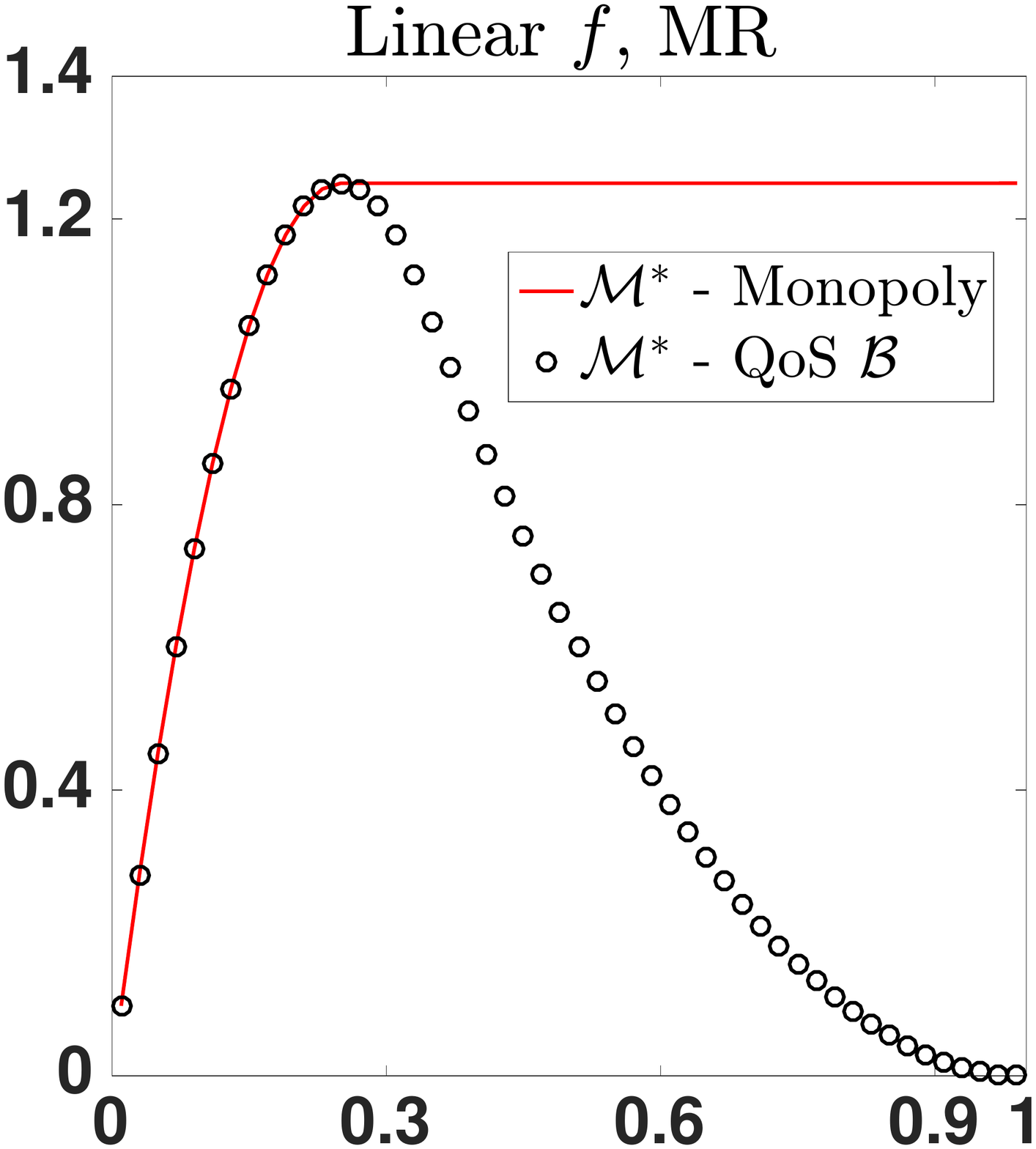}
   
        \end{minipage}
    \caption{Linear  function $f(\phi) = 1 - a\phi.$ Various equilibria or optimizers versus $e/\Lambda$ (left) and MR  at Equilibria versus $e/\Lambda$ (right), with $\Lambda = 1$, $a = 0.1$}  
    \label{fig:Lin_function}

    \centering
    \hspace{-10mm}
    \begin{minipage}{3.7cm}
     \includegraphics[width = 4.5cm, height = 4cm]{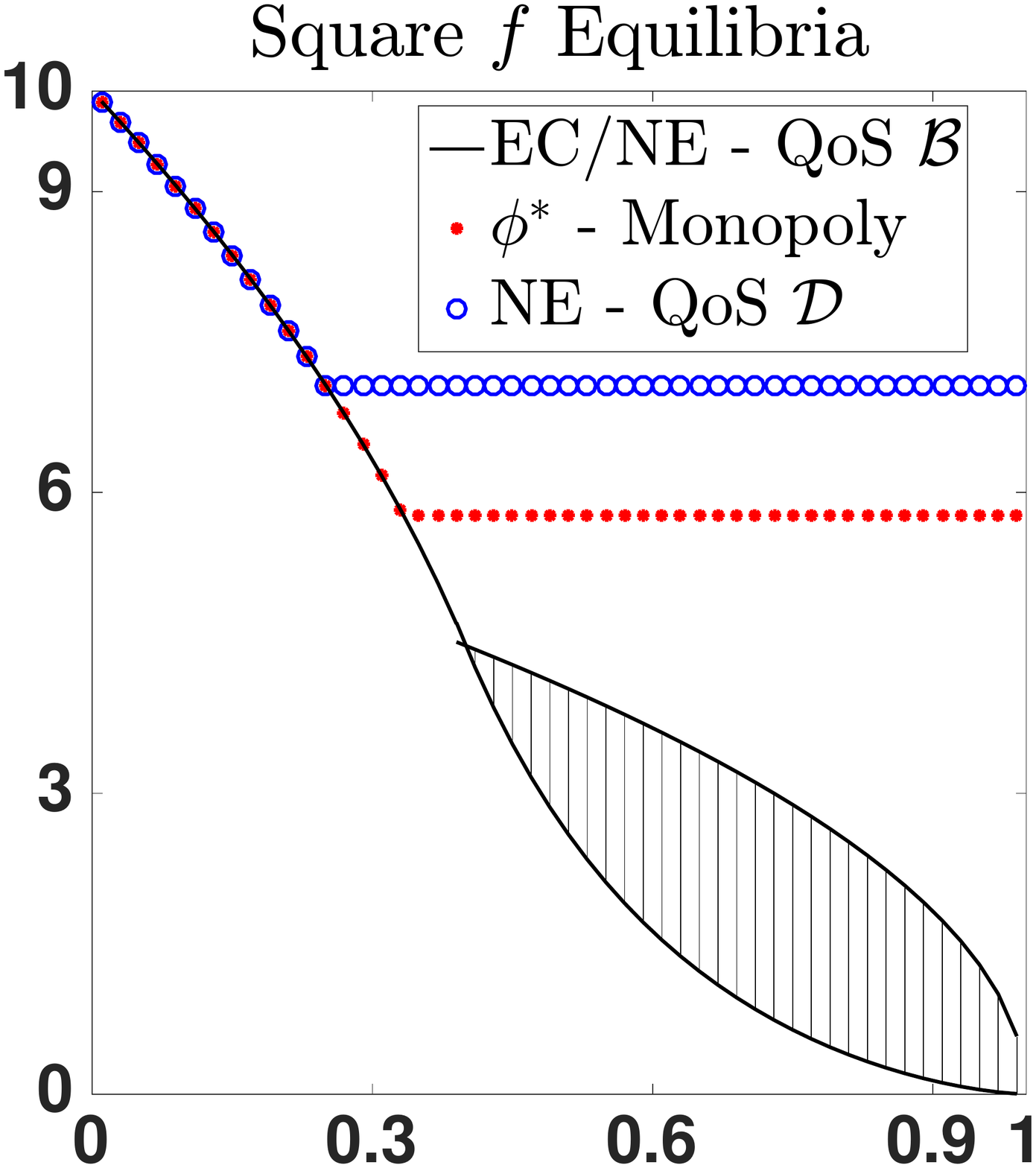}

        \end{minipage}
        \hspace{6mm}
      \begin{minipage}{3.7cm}
      \vspace{-2mm}
   \includegraphics[width = 4.5cm, height = 4cm]{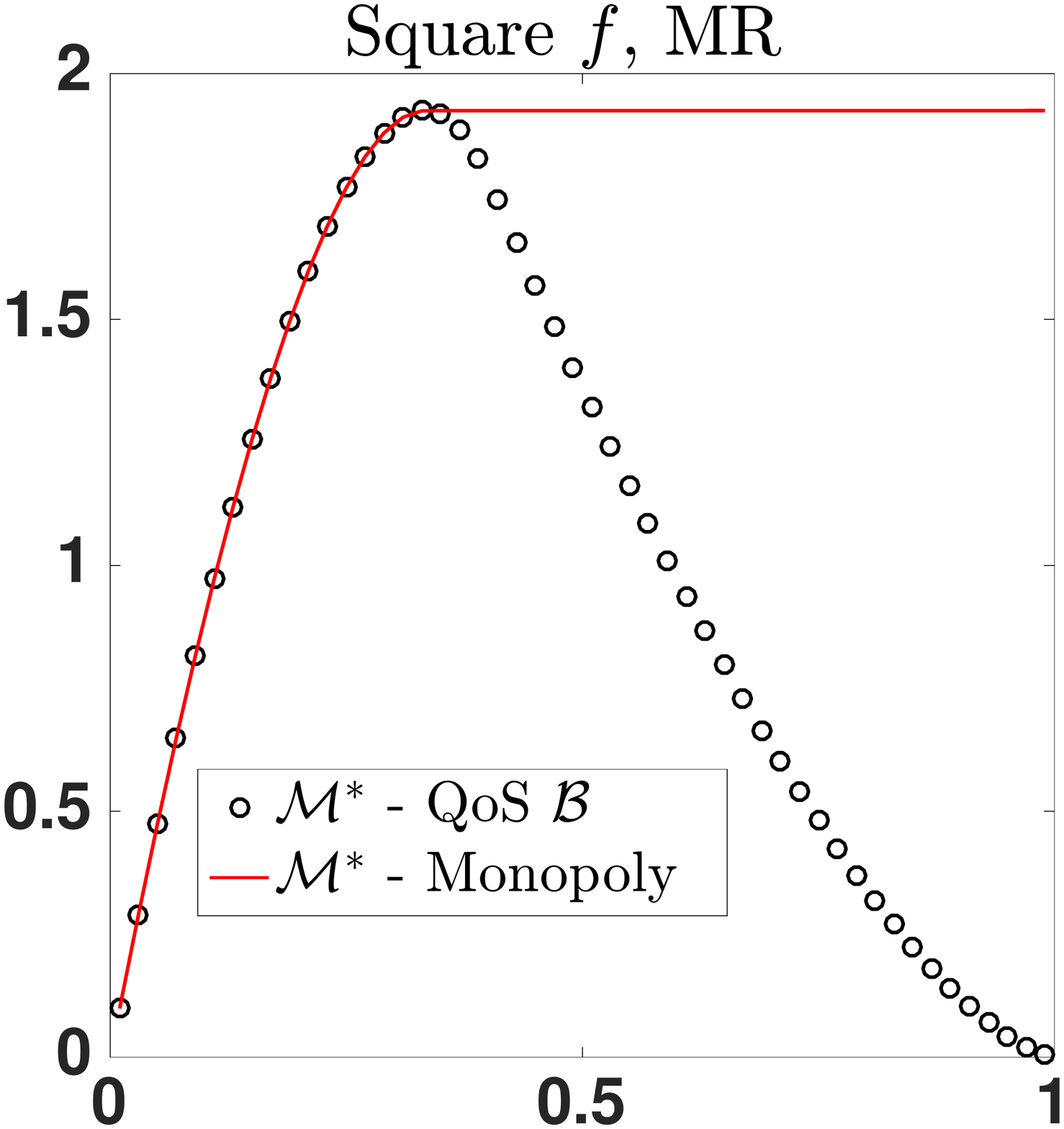}
 
        \end{minipage}
    \caption{Square  function $f(\phi) = 1 - (a\phi)^2.$ Various equilibria or optimizers versus $e/\Lambda$ (left) and MR  at Equilibria versus $e/\Lambda$ (right), with $\Lambda = 1$, $a = 0.1$}
    \label{fig:Square}
    \vspace{-4mm}
\end{figure}

\section{Conclusion}
\label{sec:conclusion}
We consider competing ride hailing platforms with impatient, price sensitive passengers, and impatient and revisiting drivers. The platform on one hand strives to meet the QoS goals of its passengers so as to capture a larger market share (the total market share across platforms is conserved). On the other hand, they also aim to increase their prices, so as to  maximize their long run revenue rate. 

We analyze scenarios with competition, cooperation and monopoly in this paper and derive closed form expressions for relevant equilibria/optimizers; the accuracy of these expressions improves, as the `impatience level' of the drivers decreases. Some important findings of our study are: a) the equilibrium prices depend only upon the ratio of the passenger and driver arrival rates and the passenger price-sensitivity function (this is in line with the observation that many results in queuing systems depend just upon load factor); b) the same dependency also holds for the normalized matching revenue (per unit passenger rate) derived at various equilibria;  c) when the drivers are willing to wait for reasonably longer time, the platforms will not benefit from cooperation unless there is a competition over the market share; d) when the passengers are sensitive to the offered prices as well as the availability of drivers, there are scenarios when a Nash equilibrium does not exist. In this case, we defined a new concept called an equilibrium cycle and established the existence of the same. 

This paper motivates future work in several directions: multiple zones, competition in the presence of dynamic pricing, drivers and/or passengers who opportunistically hop between platforms, etc. We believe the BCMP-style modeling approach we have adopted is amenable to these extensions.

\vspace{12pt}

\section*{Appendix A}

\textbf{Proof of Lemma \ref{lem_exist_unique_WE}: }
We begin with the proof of the second part (monotonicity property of QoS functions) and with $Q_i = \DA_i$. 
W.l.o.g. consider $Q_1$ and any $\lambda_1 < \lambda_1'$. 
Observe, $$\prod_{a=1}^{n_1}(\lambda_1' f(\phi_1) + a \beta_1) - \prod_{a=1}^{n_1}(\lambda_1 f(\phi_1) + a \beta_1) > 0, \text{ for any } n_1>0,$$
and hence (since both the series are convergent for any $\beta_1 > 0$),
\begin{eqnarray}
\sum_{n_1=0}^{\infty} \frac{e_1^{n_1}}{\prod_{a=1}^{n_1}(\lambda_1 f(\phi_1) + a \beta_1)} -
\sum_{n_1=0}^{\infty} \frac{e_1^{n_1}}{\prod_{a=1}^{n_1}(\lambda_1' f(\phi_1) + a \beta_1)} > 0. \nonumber
\end{eqnarray}  
Thus from \eqref{eqn_pi0_expression}, $\DA_1$  is strictly increasing and continuous. 

\noindent
Now suppose $Q_i = \PB_i$. From \eqref{Eq_overall_exact} with static price policies,
\begin{align*}
\PB_i &= \DA_i + (1 - f(\phi_i)) (1 - \DA_i)
= (1 - f(\phi_i)) + f(\phi_i)\DA_i.
\end{align*}
Clearly, $\PB_i$ is an affine transformation of $\DA_i$ with positive coefficients. Therefore, $\PB_i$ is also strictly increasing and continuous. Thus  both the QoS satisfy {\bf A}.0, when $\beta_1 > 0, \beta_2 >0$.

%
Next we show the existence of unique WE for any QoS which satisfy {\bf A}.0. Define $$g(\lambda) := Q_1(\lambda) - Q_2(\Lambda-\lambda) \text{ where } \lambda \in [0,\Lambda].$$ 
By \textbf{A.}0, $g$ is a continuous function. 
Further:
\begin{enumerate}[(i)]
    \item If $g(0) < 0$ and $g(\Lambda) \geq 0$   then using Intermediate value theorem, there exists a $\lambda_1 = \lambda^* \in (0, \Lambda)$ such that $g(\lambda^*) = 0$.  Then  $(\lambda^*, \Lambda-\lambda^*)$ is WE. Further by strict monotonicity, we have uniqueness. 
    
    \item If $g(0) < 0 $ and $g(\Lambda) < 0$. By {\bf A}.0, and definition of $g$, we have $g(\lambda) < 0$ for all $\lambda$. Hence $( \Lambda, 0)$ is the unique WE. 
    In  a similar way, when $g(0) > 0$ and $g(\Lambda) > 0$, we have that $(0, \Lambda)$ is the unique WE.

\end{enumerate}
 Rest of the cases (example $g(0) > 0$ and $g(\Lambda) \leq 0$ etc.,) follow using similar arguments.\eop
 

\textbf{Proof of Lemma \ref{lem_mr_derivation}:}
Consider a renewal process with renewal epochs being the points where the Markov process $\{Z_{t}\}_{t \ge 0}$ with $Z_t = (N_t, R_t)$, visits some state $s = (n, r)$ s.t. $n > 0$. Let the overall transition rate from state $s$ be given by $q(s)$. Then the expected length of the renewal cycle $\E(\tau(s))$ is (~\cite{balter}),
 \begin{equation}
 \label{eq_tau}
 \E[\tau(s)] = \frac{1}{q(s) \pi(s)}.
 \end{equation}
 In any renewal cycle, the platform obtains the revenue if a driver is available and the arriving passenger accepts the ride (i.e., offered price) when the system is in state $s$. Let $R(s,t)$ denote the revenue generated till time $t$, obtained while in the state $s$. This component of the reward can be obtained using the renewal process mentioned above. Towards this, let $R(s)$ represent the revenue generated in one renewal cycle of the corresponding renewal process.
 Thus, by well-known Renewal Reward Theorem (~\cite{balter}) the long-run revenue rate while in state $s$ is,
\begin{equation}
\label{eq_long_run_rev}
\lim\limits_{t \rightarrow \infty} \frac{R(s, t)}{t} \xrightarrow{\text{RRT}} \frac{\E[R(s)]}{\E[\tau(s)]} \ \text{a.s.}
\end{equation}
Observe that the expected reward generated in one renewal cycle, when in state $s$,  is the (state-dependent) price  offered by the platform ($\phi(s)$) multiplied with the probability of a passenger arriving to the system ($\nicefrac{\lambda}{q(s)}$)  and accepting the offered price  ($f(\phi(s))$), i.e., 
$$\E[R(s) ] = \phi(s) f(\phi(s)) \frac{\lambda}{q(s)}.$$

Thus, \eqref{eq_long_run_rev} can be simplified using \eqref{eq_tau}:
\begin{align}
\lim\limits_{t \rightarrow \infty} \frac{R(s, t)}{t} 
&= \frac{\phi(s) f(\phi(s)) \frac{\lambda}{q(s)}}{\frac{1}{q(s)\pi(s)}},  \nonumber \\
 &= \lambda f(\phi(s)) \phi(s)  \pi(s)  \text{ a.s.} 
 \label{eqn_revenue_using_rrt}
\end{align}

Now, define $R(t)$ as the revenue generated by the platform till time $t$. Then clearly, 
$$R(t) = \sum_{s: n>0} R(s,t).$$ 
Hence, the long-run revenue rate for the system is given by,
\begin{align*}
\lim_{t \rightarrow \infty} \frac{R(t)}{t} &= \lim_{t \rightarrow \infty} \sum_{s: n>0}\left(\frac{ R(s,t)}{t}\right),\\
& \stackrel{(a)}{=}  \sum_{s: n>0} \lim_{t \rightarrow \infty} \left(\frac{ R(s,t)}{t}\right) \stackrel{(b)}{=} \sum_{s: n>0}  \lambda f(\phi(s)) \phi(s)  \pi(s) \text{ a.s.}
\end{align*}
where the equalities $(a)$ and $(b)$ follows from Lemma \ref{lem_exchange_of_limits_mr} and \eqref{eqn_revenue_using_rrt} respectively. \eop

 \textbf{Proof of Lemma \ref{lem_approx_MR}: } 
 \Proofs{
For this lemma, we explicitly show the dependency of some functions on $\beta$. 
For any static price policy $\phi_i$, from Lemma \ref{lem_mr_derivation} we have:
\begin{eqnarray}
\label{Eqn_MR_static}
\MR_i (\phi_i; \beta)  =    
  \lambda_i f(\phi_i) \phi_i \sum_{s_i: \niw \neq 0} {\pi}_i(s_i) 
 = 
\lambda_i f(\phi_i) \phi_i (1-\DA_i ).
\end{eqnarray}
Hence it suffices to study the limit of $\DA_i$ as $\beta \to 0$. Observe here that the summation in \eqref{eqn_pi0_expression}  is always convergent for $\beta >0.$ Pick $\bar{\beta}$ such that $(\lambda_i-\epsilon) \le \lambda_i(\beta) \le (\lambda_i+\epsilon)$ for all $\beta \le \bar{\beta}$ and for some $\epsilon>0.$

\textbf{Case 1: When $e < \lambda_i f(\phi_i)$:} Pick $\epsilon$ in the above such that $e<(\lambda_i-\epsilon) f(\phi_i)$. 
From \eqref{eqn_pi0_expression} every term in summation defining $\DA_i$ can be upper bounded for all such $\beta$ by 

\vspace{-4mm}

{\small $$
\frac{e^{n}}{\prod_{a=1}^{n}(\lambda_i(\beta) f(\phi_i) + a \beta)} \le \frac{e^{n}}{\prod_{a=1}^{n}((\lambda_i-\epsilon) f(\phi_i) )}
= \left (\frac{e}{(\lambda_i-\epsilon) f(\phi_i)} \right )^n.$$}
Thus, one can uniformly upper bound
$$
\sum_{n=0}^\infty \frac{e^{n}}{\prod_{a=1}^{n}(\lambda_i(\beta) f(\phi_i) + a \beta)} 
\le \frac{(\lambda_i-\epsilon) f(\phi_i)}{(\lambda_i-\epsilon) f(\phi_i) - e}.
$$
    Clearly the upper bounding term is summable and uniformly bounds the left hand series for any $\beta \le \bar{\beta}$. 
Hence  by convergence of each term of the series and  Dominated Convergence Theorem (DCT)  \cite{protter}, 
\begin{equation}
  \DA_i (\phi_i; \lambda_i(\beta),\beta) \to 1 - \frac{e}{\lambda_i f(\phi_i)} := \DAl_i  (\phi_i),  \mbox{ as, } \beta \to 0. \label{Eqn_DA_i_conv_one}
\end{equation}

\textbf{Case 2: When $e \geq \lambda_i f(\phi_i) > (\lambda_i+\epsilon) f(\phi_i)$,} the following series   diverges:
\begin{align*}
& \sum_{n=0}^\infty \frac{e^{n}}{\prod_{a=1}^{n}((\lambda_i+\epsilon) f(\phi_i)) }  
= \sum_{n=0}^\infty  \left (\frac{e} {(\lambda_i+\epsilon) f(\phi_i)} \right  )^n 
\end{align*}
Thus, 
for every $\delta > 0$, there exists an $N_\delta$   such that 
\begin{align*}
& \sum_{a=0}^{N_\delta}
\left ( \frac{e }{( \lambda_i + \epsilon) f(\phi_i)  }  \right )^a \ge \frac{1}{\delta}
\end{align*}
For this $N_\delta$
choose ${\bar \beta}_\delta$ small enough such that by continuity  (of finite sum) for all $\beta \le {\bar \beta}_\delta$ 
\begin{align*}
  \sum_{n=0}^{N_\delta}\left( \frac{e^{n}}{\prod_{a=1}^{n}((\lambda_i + \epsilon) f(\phi_i) + a \beta)}\right) \ge \frac{1}{2\delta} 
\end{align*}
which by positivity of term implies (for all $\beta \le \bar{\beta}$)

\vspace{-4mm}

{\small \begin{align*}
  \sum_{n=0}^{\infty}\left( \frac{e^{n}}{\prod_{a=1}^{n}(\lambda_i(\beta) f(\phi_i) + a \beta)}\right) &\ge \sum_{n=0}^{\infty}\left( \frac{e^{n}}{\prod_{a=1}^{n}((\lambda_i+\epsilon) f(\phi_i) + a \beta)}\right)  \\
  &\ge \frac{1}{2\delta}   \mbox{, and hence, } 
 \DA_i   \le  2\delta .
\end{align*}}
Thus $\DA_i  (\phi_i; \lambda_i(\beta),\beta) \to 0 := \DAl_i  (\phi_i)$ as $\beta \to 0 $. Hence and further using \eqref{Eqn_MR_static} and \eqref{Eqn_DA_i_conv_one},   the result follows. \eop
}
{
This proof follows using standard analysis based arguments and is omitted due to space constraints (available in \cite{TR}). \eop
}

\textbf{Proof of Theorem \ref{thm_monopoly_opt}: }
We begin with the proof of the first part by showing convergence along a sub-sequence to the unique optimal price of the limit system, while the uniqueness of the optimal price for the limit system is established later.


\noindent
\textbf{Convergence of Optimal Prices:}  The matching revenue of the system for any $(\beta, \phi)$, where $\phi$ is a static price policy, by Lemma \ref{lem_mr_derivation} and definition of $\DA$ as in  \eqref{eqn_pi0_expression} equals:
$$
\MR(\phi;\beta)   \ = \   \frac{\Lambda}{2}f(\phi) \phi(1-\DA(\phi; \beta)).
$$
From \eqref{eqn_pi0_expression},
it is easy to observe that the above function is jointly continuous in $\phi$ and $\beta$.  
The optimal price $\phi^*$
for any $\beta \in (0, \infty)$ maximizes the above   function, i.e.,
$$
\phi^* \in 
 \Phi^*(\beta) := \arg \max_{\phi \in [0,\phi_h]} \MR(\phi;\beta).
$$
Clearly, the domain of optimization $[0,\phi_h]$ is  compact and is the same for all $\beta$ and hence, the hypothesis of Maximum Theorem \cite[Theorem 9.14]{maximum} is satisfied. Thus the correspondence $\beta \mapsto \Phi^*(\beta)$ is compact and upper semi-continuous. 

Consider any sequence $\beta_n \to 0$.
 Consider one optimizer, from $\Phi^*(\beta)$ for each $\beta_n$, and call it $\phi_n^*$. By \cite[Proposition 9.8]{maximum}, there exists  
 a sub-sequence of $\{\phi^*_n\}$,  which converges to   the unique optimal price of the limit system.

\noindent
\textbf{Unique Optimal Price   at limit:}
By Lemma \ref{eq_approx_MR}, the  matching revenue of  the limit system equals (recall $\phiunderbar = f^{-1}\left(\nicefrac{2e}{\Lambda}\right)$),
\begin{eqnarray} \label{eqn_monopoly_mr}
\MRl(\phi) &=&  
e \phi  \mathds{1}_{ \{ \phi \in  [ 0, \phiunderbar  ) \} } + M(\phi)  \mathds{1}_{ \{   \phi \in [\phiunderbar, \phi_h ] \} }  \\
\mbox { with } M(\phi) &: =& \frac{\Lambda}{2} \P(\phi). \nonumber
\end{eqnarray}
Observe that $M(\phiunderbar) = e \phiunderbar $   and hence $\MRl$ is continuous over $\phi \in [0, \phi_h]$. Further, from Lemma \ref{lem_monopoly},  the function $M  $ is strictly concave and let $\phifp := \arg \max_\phi M(\phi)$ be its unique maximizer. Thus i) if $\phifp \leq \phiunderbar$ 
then the unique maximizer of $\MRl$ is $\phiunderbar$ ($M$ is decreasing for  $\phi \ge \phiunderbar$); else 
ii) 
when $\phifp > \phiunderbar$  the unique maximizer of $\MRl$ coincides with that of $M$, 
$
\phifp
$ (see     \eqref{eqn_monopoly_mr}). \eop

\textbf{Proof of Lemma \ref{lem_pio_WE}}: 
From Lemma \ref{lem_exist_unique_WE}, we know that there exist an unique WE with $\DA_i$ as a WE metric. Our aim is to find the passenger arrival rate split  at WE $(\lambda_1, \lambda_2)$ s.t . it minimizes the function $(\DA_1 - \DA_2)^2$. Suppose if we assume \begin{eqnarray} \label{eqn_WE_Split_eqn}
\left(\lamo,\lamt\right)  = \left(\frac{\Lambda f(\phi_2)}{f(\phi_1) + f(\phi_2)}, \frac{\Lambda f(\phi_1)}{f(\phi_1) + f(\phi_2)}\right).
\end{eqnarray}
Then $(\DA_1 - \DA_2)^2 = 0$ and we know that the minimum of $(\DA_1 - \DA_2)^2$ cannot be lesser than zero. Therefore the unique minimizer of the function is as given in \eqref{eqn_WE_Split_eqn}. \eop

\textbf{Proof of Theorem \ref{thm_pio_NES}: } From \eqref{eq_approx_MR}, the limit MR of platform when it chooses strategy $\phi_i$ against any $\phi_{-i}$ equals,
 \begin{eqnarray}
 \label{eq_MR_thm}
    \MRl_i(\VPhi) = \min \{ e\phi, h(\VPhi) \} \text{ where } h(\VPhi) = \frac{\Lambda f(\phi_i)f(\phi_{-i})\phi_i}{f(\phi_i)+f(\phi_{-i})}.
\end{eqnarray}
Next we define the price at which the above two functions become equal (if at all),

\vspace{-4mm}

{\small \begin{equation}
    \hat{\phi}(\phi_{-i}) = \begin{cases}
    \phi_h & \text{ if } h(\phi_i,\phi_{-i})>e\phi_i \text{ for all } \phi_i, \\
    0 & \text{ if } h(\phi_i,\phi_{-i})<e\phi_i \text{ for all } \phi_i \\
    \inf_{\phi_i} \left \{ h(\phi_i,\phi_{-i})=e\phi_i\right \} & \text{ else.}
    \end{cases}
\end{equation}}
From \eqref{eq_MR_thm},  $\MRl_i (\phi_i, \phi_{-i}) = h (\phi_i, \phi_{-i})$ for all $\phi_i \ge \hat{\phi}(\phi_{-i})$, and this is  because   $\phi_i \mapsto h(\phi, \phi_{-i})/\phi_i$ is strictly decreasing (the first derivative is always negative). 

By simple algebra,  one can show that when  $\phi_{-i} = \phiunderbar$, we have $\hat{\phi} (\phi_{-i}) =\phi_{-i} $ (and this is the only such $\phi_{-i}$).  Thus we have the proof of part (i), if $\phi_i \mapsto h (\phi_i, \phiunderbar) $ is decreasing for all $\phi_i \ge \phiunderbar$. 
Towards this, consider the relevant partial derivative,
\begin{eqnarray}
n(\phi_{-i})  =  \left . \frac{\partial h}{\partial \phi_i} \right |_{\hat{\phi},\phi_{-i}}   \
 =  \ \left(\frac{e}{\Lambda} \right)^2 \Lambda \hat{\phi} \left( \frac{\Lambda}{e} + \frac{\hat{\phi}f'(\hat{\phi})}{(f(\hat{\phi}))^2} \right), \hspace{3mm}
\end{eqnarray}
where the last equality follows using 

\vspace{-4mm}
{\small 
\begin{eqnarray}
\label{eq_derivative_h}
\frac{\partial h}{\partial \phi_i} = \Lambda f(\phi_{-i}) \phi_i \left[ \frac{(f(\phi_i))^2+f(\phi_i)f(\phi_{-i})+f(\phi_{-i})f'(\phi_i)\phi_i}{(f(\phi_i)+f(\phi_{-i}))^2}\right],\\
    \mbox{ and, } \frac{ f(\hat{\phi})f(\phi_{-i})}{f(\hat{\phi})+f(\phi_{-i})} = \frac{e}{\Lambda},
 \mbox{ 
which implies }
 f(\phi_{-i}) = \frac{e f(\hat{\phi})}{\Lambda f(\hat{\phi})-e}. \nonumber   
\end{eqnarray}}
Hence part (i) follows,  because: i) $\hat{\phi}(\phiunderbar) = \phiunderbar$;
ii) $n(\phiunderbar) = \Lambda\phiunderbar d(\phiunderbar)/4 \le 0 $; and
iii) hence and further by Lemma \ref{lem_negative_partial_derivative}, $h$ is decreasing  beyond $\phiunderbar$ (since $\hat{\phi}(\phiunderbar) = \phiunderbar$ which implies $h(\phiunderbar,\phiunderbar)=e\phiunderbar$ and beyond $\phiunderbar$, $\MRl_i = h$); and finally iv) BR against $\phiunderbar$ is $\phiunderbar$.





\textbf{Parts (ii) and (iii):} Next consider the case when $d(\phiunderbar)>0$. 
This implies $n(\phi_{-i}) > 0$ for all $\phi_{-i} \ge \phiunderbar$, as 
$$
\frac{n (\phi_{-i}) f({\hat \phi})^2}{e^2/\Lambda \hat{\phi} } = \frac{\Lambda}{e} f({\hat \phi})^2 + \hat{\phi} f'({\hat \phi})
$$increases with $\phi_{-i}$   and coincides exactly with $d(\phiunderbar)$ at   $\phi_{-i} = \phiunderbar$; the first statement is true because $\phi_{-i} \mapsto {\hat \phi}(\phi_{-i})$  is a decreasing function (derivative always negative) and hence $\phi_{-i} \mapsto f({\hat \phi})$ is  increasing (as in Lemma \ref{lem_monotonicity_of_h_1_function}). 

In summary  $n(\phi_{-i}) > 0$ for all $\phi_{-i} \ge \phiunderbar$. This implies that $h $ and hence $\MRl_i$ is increasing at least initially after $\phi_{i} \ge \hat{\phi}(\phi_{-i})$. 
Thus we have a symmetric NE if there exists a $\phi_{-i} > \phiunderbar$ for which 
$$
\left . 
\frac{\partial h}{\partial \phi_i}  \right |_{   (\phi_{-i}, \phi_{-i})} =   0
$$
From \eqref{eq_derivative_h},
  the above partial derivative at $(\phi_{-i}, \phi_{-i})$
equals $
l(\phi_{-i}) d(\phi_{-i})$, for some positive  function $l (\phi_{-i})=\Lambda \phi_{-i}/4$. Thus the zero $\phido$ of $d$ forms a part of  the symmetric NE (when it exists), as further the relevant second
derivative can be shown to be negative. 
If there are no such  zeros clearly $(\phi_h, \phi_h)$ is the NE 
as then the partial derivative $\partial h/ \partial \phi |_{\phi_{-i} = \phi_h}$ is always positive. \eop


\textbf{Proof of Theorem \ref{thm_WE_MR_PB}:} We define the WE and MR of the limit system uniquely such that we have continuity as $\beta \to 0$. 
The procedure is as follows. We first find the minimizers of the objective function in \eqref{eqn_existence_uniq_WE_main} by treating $\PBl$ of \eqref{Eq_overall_aproxx} as limit QoS. 
Say at a particular $\VPhi$,  these minimizers are unique. 
Then  continuity w.r.t. $\beta$ follows from Maximum theorem \cite{maximum} as in the proof of Theorem \ref{thm_monopoly_opt}, with additional advantage due to unique optimizers: the upper semi-continuity of the correspondences of optimizers implies  continuity when they are further unique (~\cite{maximum}). If the minimizers of limit system are  not unique, we pick an appropriate one from among, so as to ensure continuity (such a choice is again possible, because by Lemma \ref{lem_exist_unique_WE}, WE are unique for $\beta > 0$ and hence   $\lim_{\beta \to 0} W(\beta)$ always exists). We call optimizers of $\PBl$ as limit WE.
Once limit WE are found, we obtain limit MR using Lemma \ref{lem_approx_MR}, and show it equals the ones provided in Table \ref{tab:WE_MR_PB}. 
We proceed with the proof in the following steps.  

\textbf{Case 1: When $\phi_1=\phi_2$:} From \eqref{Eq_overall_aproxx}, it is easy to see that $\lamol = \Lambda/2$
is a limit WE. We prove the uniqueness by contradiction. If possible, w.l.o.g. say there exists another limit WE with $\epsilon>0$,
$$
(\lambda'_1(\VPhi),\lambda'_2(\VPhi)) = \left( \frac{\Lambda}{2}-\epsilon,\frac{\Lambda}{2}+\epsilon \right).
 $$
  One can easily verify from \eqref{Eq_overall_aproxx} that the blocking probabilities of the platform $1$ and $2$ are $1-f(\phi_1)$ and $1-\nicefrac{e}{\lambda'_2(\VPhi)}$ respectively which can never be equal  since, $$\frac{e}{\lambda'_2(\VPhi)}<f(\phi_2)=f(\phi_1).$$
So, we have  unique minimizer of \eqref{eqn_existence_uniq_WE_main} for   limit system and for such $\VPhi$ (when $\PBl$ are equal, the optimal value  of  \eqref{eqn_existence_uniq_WE_main} is $ 0$).
By Lemma \ref{lem_approx_MR}, the matching revenue of the limit system equals,
 $$
    \MRl_i(\VPhi) = e\phi_i\mathds{1}_{\phi_i \in [0,\phiunderbar)} + \frac{\Lambda}{2} f(\phi_i)\phi_i \mathds{1}_{\phi_i \in [\phiunderbar,\phi_h]}.
$$

\textbf{Case 2: When $\phi_1 >  \phi_2$:}

(a) When $\phi_1, \phi_2 \in [0, \phiunderbar)$: From the definition of $\phiunderbar$ (in \eqref{Eqn_Mphi_PhistarMPhi}), we have $\nicefrac{2e}{\Lambda f(\phi_i)} < 1$ for any $i \in \{1,2\}$. Then, from  \eqref{Eq_overall_aproxx}, it is easy to see that the limit WE again equals, 
$\lamol = \Lambda/2$. To prove the uniqueness, say there exists another WE with $\epsilon>0$,
$$
(\lambda'_1(\VPhi),\lambda'_2(\VPhi)) = \left( \frac{\Lambda}{2}-\epsilon,\frac{\Lambda}{2}+\epsilon \right).
 $$
 
 (i) When $\PBl_1'=1-f(\phi_1)$ and $\PBl_2'=1-f(\phi_2)$:
 The above blocking probabilities can never be equal  since $f(\phi_1)\neq f(\phi_2)$.
 
 (ii) When $\PBl_1'=1-\nicefrac{e}{\lambda'_1(\VPhi)}$ and $\PBl_2'=1-\nicefrac{e}{\lambda'_2(\VPhi)}$: These cannot be equal since $\lambda'_1(\VPhi)\neq \lambda'_2(\VPhi)$.
 
(iii) When $\PBl_1'=1-f(\phi_1)$ and $\PBl_2'=1-\nicefrac{e}{\lambda'_2(\VPhi)}$: This case is also not possible since,
$$\frac{e}{\lambda'_2(\VPhi)}<\frac{e}{\lamtl}=\frac{e}{\lamol}<f(\phi_1).$$
Similar logic follows for the remaining sub-case and by Maximum theorem, $ W(\beta) \to W(0) = \Lambda/2$, the unique limit WE, as given in Table \ref{tab:WE_MR_PB}.
%
By Lemma \ref{lem_approx_MR}, the matching revenue of the limit system can be set as,
$
    (\MRl_1(\VPhi),\MRl_2(\VPhi)) = (e\phi_1,e\phi_2).
$


(b) When $\phi_1 \in [\phiunderbar, \phibar)$ with $\phibar = f^{-1}\left(\frac{e}{\Lambda}\right)$, the limit WE   is given by, \begin{equation}
\label{eq_WE_b}
    \lamol   =  \left(\Lambda  - \nicefrac{e}{f(\phi_1)}\right).
\end{equation} 
 One can verify the above by directly substituting   in \eqref{Eq_overall_aproxx}; the resulting $\PBl$ are again equal and are given by $1- f(\phi_1)$ and $1 - \frac{e}{\lamtl}$ respectively. 

Next, we show that the above WE is the unique one. Towards this we first study $\PBl_1(\lambda_1)$ and $\PBl_2 (\Lambda-\lambda_1)$ of \eqref{Eq_overall_aproxx} as a function of $\lambda_1$.  With $x_1:=\nicefrac{e}{  f(\phi_1)}$ and 
$x_2:=\Lambda- \nicefrac{e}{  f(\phi_2)},$ 
\begin{align*}
    \PBl_1(\lambda_1) &= 1- \left(f(\phi_1) \mathds{1}_{ \lambda_1 \le x_1}  + \frac{e}{\lambda_1}   \mathds{1}_{ \lambda_1 > x_1}\right)  \\
     \PBl_2(\lambda_1) &=1- \left(f(\phi_2) \mathds{1}_{ \lambda_1 \ge x_2}  + \frac{e}{(\Lambda-\lambda_1)}   \mathds{1}_{ \lambda_1 < x_2}\right).
\end{align*}
Thus from \eqref{eq_WE_b} and $\phi_1 \geq \phiunderbar$, $\lamol \le x_1$ and $\lamol < x_2$.  
And thus there can be no other $\lambda_1$ at which the two $\PBl$ become equal: from the above equation, for all $\lambda_1 > \lamol$, $\PBl_1 (\lambda_1) > \PBl_2 (\lambda_1) $; the vice versa is true for $\lambda < \lamol$.

By Lemma \ref{lem_approx_MR} again, the matching revenue of the limit system can be set as,
$$
(\MR_1(\VPhi),\MR_2(\VPhi)) = ((\Lambda f(\phi_1)-e)\phi_1,e\phi_2).
$$
Observe $\lamol f(\phi_1)< e$ because $\phi_1 \ge \phiunderbar$ and so $f(\phi_1) \le 2e /\Lambda$.

(c) When $\phi_1 \in [\phibar, \phi_h]$, from \eqref{Eq_overall_aproxx} , $ \PBl_1 = 1 - f(\phi_1)$ irrespective of $\lamo$ while $\PBl_2 = \min\left\{ 1 - \frac{e}{\lamtl}, 1-f(\phi_2) \right \}$ where $\lamtl$ has to be found. Observe that for any $\lamtl$
$$
1-f(\phi_1) >1-f(\phi_2)\ge \min \left\{1 - \frac{e}{\lamtl}, 1 - f(\phi_2) \right\},
 $$
since $f(\phi_1) < f(\phi_2)$ and hence it is not possible to equate the blocking probabilities of the two platforms. 
In this case any $\lambda_1 $ optimizes \eqref{eqn_existence_uniq_WE_main}. 
Thus, we consider the following as limit WE:  
$$(\lamol, \lamtl) = \left(0, \Lambda \right),$$
and show that the selection ensures the required continuity properties. Towards this, consider any $\beta_n \to 0$ and consider corresponding $\lim_{\beta_n \to 0} W(\beta_n)$. We claim that $W(\beta_n) \to 0$, which we prove case-wise as below:
  
\textbf{When $f(\phi_2) \le \nicefrac{e}{\Lambda}$:} 
Consider any $\epsilon >0.$
By   Lemma \ref{lem_approx_MR}, with $\lambda_i(\beta) = 0$
for all $\beta$,     
  there exists a $\beta_\epsilon$ such that 
  $$\DA_1 (\phi_1; 0, \beta) \le \epsilon/f(\phi_1)  \mbox{ and hence }
  \PB_1(\phi_1; 0,\beta) < 1-f(\phi_1) + \epsilon$$
  for all  $\beta \le \beta_\epsilon$ (see   \eqref{Eq_overall_exact}). In similar lines,  
  $$\PB_2(\phi_2; \Lambda,\beta) > 1-f(\phi_2)-\epsilon, \mbox{ when } \beta \le \beta_\epsilon.$$  This implies $\PB_1(\lambda;\beta) < 1-f(\phi_1) + \epsilon$ and $\PB_2(\lambda;\beta) \ge 1-f(\phi_2)-\epsilon$ respectively, for all $\lambda$. 
For all such $0 < \beta \le \beta_\epsilon$ from \eqref{eqn_existence_uniq_WE_main},
$W(\beta)  = (0, \Lambda)$, if $\epsilon < (f(\phi_2) - f(\phi_1))/4$.

One can give similar arguments for the other case, now with $\PB_2(\phi_2;\Lambda,\beta) \to 1-\nicefrac{e}{\Lambda f(\phi_2)}$.
This completes the proof of claim. 
   Again from Lemma \ref{lem_approx_MR}, the matching revenue can be set to,
$$\hspace{6mm}
(\MRl_1(\VPhi),\MRl_2(\VPhi)) = (0,\min \{\Lambda f(\phi_2),e\} \phi_2). \hspace{6mm} \mbox{ \eop}
$$



{\bf Proof of Theorem \ref{thm_PB_NE_phiunderbar}:}
We directly find the best response (BR) against $\phi_2 = \phiunderbar$ (w.l.g.). 
From Table \ref{tab:WE_MR_PB}, the matching revenue of platform $1$   is given below

{\small\begin{equation*}
\MR_1 (\phi_1,\phiunderbar) =
e\phi_1 \mathds{1}_{   \phi_1 \in \left[0,  \phiunderbar,   \right) }
+ m(\phi_1) 
\mathds{1}_{ \phi_1 \in \left[ \phiunderbar, \phibar  \right) },  \
m(\phi) = (\Lambda f(\phi) - e )\phi.
\end{equation*}}
By   Lemma \ref{lem_monopoly}, 
and the given hypothesis ($\phimr < \phiunderbar$), we have that $m$ is strictly decreasing for $\phi \ge \phiunderbar$ and $m(\phiunderbar) = e\phiunderbar$ and hence the BR of player 1 is $\phiunderbar$.  \eop

{\bf Proof of Theorem \ref{Thm_PB_EC}: }
Consider two platforms offering price vector $\VPhi = (\phi_i,\phi_{-i})$ such that $\phi_{-i} \in \left[\eL, \phimr\right]$.
From Table \ref{tab:WE_MR_PB}, the matching revenue of platform $i$ against such $\phi_{-i}$ (observe $\phi_{-i} \le \eU \le \phibar$)   is given by (with $m(\phi) = (\Lambda f(\phi)  - e) \phi$),
\begin{equation} \label{eqn_MR_equilibrium_cycle_proof}
\MR_i(\phi_i, \phi_{-i}) =
 \begin{cases}
 e \phi_i & \text{ if } \phi_i < \phi_{-i} \\
 \frac{\Lambda}{2} f(\phi_{i}) (\phi_{i}) & \text{ if }  \phi_i = \phi_{-i} \\
 \max\{m(\phi_i), 0\} & \text{ if }  \phi_i > \phi_{-i}.
\end{cases}   
\end{equation}
Let $\mstar := m(\phimr).$ Also observe by definition that $\mstar = e \eL$

\textbf{Part (i): }
We begin with the proof of the first part of Definition \ref{def_equi_cycle} of the equilibrium cycle. Begin with 
a $\phi_i  < \eL$, and then  $\MR_i(\phi_i, \phi_{-i}) = e \phi_i$ from   \eqref{eqn_MR_equilibrium_cycle_proof}. We will find the required $\phi_i'$ depending upon $\phi_{-i}$.

\begin{enumerate}
    \item If $\phi_{-i} = \eL$,  
    then with   $\phi_i' = \phimr$ we have $$
     \MR_i(\phi_i', \phi_{-i})  = \mstar = e \eL > e \phi_i = \MR_i(\phi_i, \phi_{-i}) .
    $$

    \item If $\phi_{-i} \in \left( \eL, \phimr \right]$, set  $\phi_i'= \eL$. Then from   \eqref{eqn_MR_equilibrium_cycle_proof},
    $$\MR_i(\phi_i', \phi_{-i}) = e \phi_i' >  e \phi_i = \MR_i(\phi_i, \phi_{-i}) .$$

    
\end{enumerate}
Now consider that $\phi_i > \eU$ and set $\phi_i'= \eU$. Clearly from \eqref{eqn_MR_equilibrium_cycle_proof},
%
$$\MR_i(\phi_i', \phi_{-i}) = \begin{cases}
\mstar & \text{if } \phi_{-i} < \phi_i' \\
\frac{\Lambda}{2}f(\phimr) \phimr & \text{if } \phi_{-i} = \phi_i'
\end{cases}$$
From Lemma \ref{lemma_oscillatory},  $\MR_i(\phi_i', \phi_{-i}) \geq \mstar$, and
since $\phi_i > \phi_i'$, we have
$$
\MR_i(\phi_i', \phi_{-i}) \geq \mstar > max\{0, m(\phi_i)\} = \MR_i(\phi_i, \phi_{-i}).
$$
This proves the first part of the definition. 


\textbf{Part (ii): } Now consider any $\phi_1,\phi_2 \in [\eL, \eU]$. By symmetry of limit system, sufficient
to consider the following   cases

\begin{enumerate}
    \item  If $\phi_1 > \phi_2$, pick any
  $\phi_2' \in (\phi_2, \phi_1)$.  Then 
  $$
  \MR_2(\phi_2', \phi_1 ) = e\phi_2' >  e\phi_2 = \MR_2(\phi_2, \phi_1). 
$$  
 
 \item  If $\phi_1 = \phi_2   \neq \eL$, then\footnote{ As $\eU > \phiunderbar$, by definitions,  $e \eL = m(\eU) > m(\phiunderbar) = e \phiunderbar$, thus $\eL > \phiunderbar $, and thus by Lemma \ref{lemma_oscillatory},  
   $\Lambda f(\phi_2) \phi_2  /2  < e\phi_2 $.
   Thus, $  \max \left \{\frac{\Lambda f(\phi_2)\phi_2}{2e}, \eL \right \}    < \phi_2$.
   } let $\phi_2' \in \left( \max \left \{\frac{\Lambda f(\phi_2)\phi_2}{2e}, \eL \right \} , \phi_2\right) \subset [\eL, \phi_2)$. Then
 $$
 \MR_2(\phi_2', \phi_{1}) = e\phi_2' > \frac{\Lambda f(\phi_2)\phi_2}{2} =  \MR_2(\phi_2, \phi_{1}).
 $$

 \item 
 If $\phi_1 = \phi_2 = \eL$, let $\phi_2' = \phimr$. By Lemma \ref{lemma_oscillatory} (as $\eL > \phiunderbar$),
 $$
 \MR_2(\phi_2', \phi_1) = \mstar = e \eL 
 >\frac{\Lambda}{2} f\left(\eL\right) \eL = \MR_2(\phi_2, \phi_1). 
 $$
 \end{enumerate}
 This completes the proof of part2 of Definition \ref{def_equi_cycle}. \eop

\textbf{Proof of Theorem \ref{thm_epsilon_NE}}
If $e \ge \Lambda$  then for all $\phi$, $ \frac{e}{\Lambda} \ge 1 \geq f(\phi)$ which implies that $\phi \in [\phibar,\phi_h]$. Therefore, for any $\phi$, $ \frac{e}{\Lambda} \ge f(\phi)$. In particular consider $ \delta $ given in hypothesis, and we will show that $\delta$ is in $\epsilon$-BR. 

From Table \ref{tab:WE_MR_PB}, the matching revenue against this $\delta $ is:

\begin{equation*}
\MR_{\delta}(\phi) =
 \begin{cases}
 0 & \text{ if } \phi > \delta \\
 \frac{\Lambda}{2} f(\delta) \delta & \text{ if }  \phi = \delta \\
 \Lambda f(\phi) \phi & \text{ if }  \phi < \delta
\end{cases}   
\end{equation*}
Thus the result.  (Here note that, $\sup_{\phi \leq \delta} \Lambda f(\phi) \phi < \epsilon$ then $\frac{\Lambda}{2} f(\delta) \delta < \epsilon$. Therefore the player is $\epsilon$-indifferent between these two payoffs.)  \eop

\ignore{\begin{thm}
\label{thm_new}
If $\phifp, \phimr, \phido \in (0, \phi_h)$ then $\phifp <  \phimr < \phido$.
\end{thm}

\textbf{Proof of Theorem \ref{thm_new}: }We know that $\phimr = \arg \max_{\phi} \Lambda f(\phi) \phi - e\phi$, $\phifp = \arg \max_{\phi} \Lambda f(\phi) \phi$ and $\phido$ is the solution to the fixed point equation $2f(\phi) + \phi f'(\phi) = 0$. Now, $m'(\phi) = \Lambda (f(\phi) + f'(\phi) \phi ) - e$ and $M'(\phi)= \Lambda (f(\phi) + f'(\phi) \phi ) $ We will prove the theorem in two steps:

step 1 (Claim: $\phifp < \phimr$)

Observe that $ M'(\phi) = m'(\phi) + e$. The function $f(\phi)$ is strictly decreasing. Since $f'(\phi) < 0$, $\phi f'(\phi)$ is also a decreasing function of $\phi$. Therefore $M'(\phi) = \Lambda (2f(\phi) + \phi f'(\phi))$ and $m'(\phi)$ are also a strictly decreasing functions of $\phi$. }

\section*{Appendix B}

\begin{lemma} \label{lem_exchange_of_limits_mr}
The long-run revenue rate of the system is equivalent to the long-run revenue rate of the system in state $s$ when summed over all possible states, i.e.,
$$\lim_{t \rightarrow \infty} \sum_s \frac{R(s,t)}{t} =  \sum_s \lim_{t \rightarrow \infty} \frac{R(s,t)}{t} \ \ \text{a.s.}$$
\end{lemma}
\noindent
\textbf{Proof: }Consider the following:
\begin{equation}
\label{eq_epsilon}
 \epsilon_t :=  \left| \sum_s \frac{R(s,t)}{t} - \sum_s R^*(s) \right|, \text{ where } R^*(s) = \lim_{t \rightarrow \infty} \frac{R(s,t)}{t}. 
\end{equation} 
To prove the required result, one needs to show that $\epsilon_t \to 0$ as $t \to \infty$ a.s.
Next, we define the following:
\begin{align*}
C_b &:= \{s: n+r \ge b\} \text{ for every }b \in \mathds{Z}^+ \text{ and,} \\
\mathcal{A}_b &:= \left \{\omega : \frac{V_t(C_b)(\omega)}{t} \rightarrow \pi(C_b) \right \} \cap \left \{\frac{R(s,t)}{t} \rightarrow R^*(s)  \ \forall \ s \right \},
\end{align*}
where $V_t(C_b)$ is the number of visits to the set $C_b$ till time $t$. By Law of Large Numbers (LLN) in standard textbooks (for example, \cite[Theorem 1.10.2]{norris}), we know that as $t \to \infty$, time-average ($\nicefrac{V_t(C_b)}{t}$) converges to stationary measure  ($\pi(C_b)$) almost surely  and hence the first set has probability $1$. Further, by RRT (as shown in Lemma \ref{lem_mr_derivation}), $\Prob(\mathcal{A}_b) = 1 $ for all $b$.

Let $\mathcal{A} = \cap_b \mathcal{A}_b$. In other words, we have
\begin{equation*}
\mathcal{A} := \left \{\omega : \frac{V_t(C_b)(\omega)}{t} \rightarrow \pi(C_b) \ \forall \ b \right \} \cap \left \{\frac{R(s,t)}{t} \rightarrow R^*(s)  \ \forall \ s \right \}.
\end{equation*}

Using continuity of probability \cite[Theorem 2.4]{protter}, we get that $\Prob(\mathcal{A}) = 1$.

Now, the final step is to show that for all $\tilde{\delta} > 0$, there exists a $\tilde{T}_\delta$ such that $\epsilon_t \le \tilde{\delta}$ when $t \ge \tilde{T}_\delta$. Towards this, consider any $\omega \in \mathcal{A}$ and for any $\delta$, choose a $b_\delta >0$ such that
$$
\sum_{s: n+r > b_\delta} \pi(s) < \frac{\delta}{4}.
$$
This is possible because the stationary distribution exists (see Lemma \ref{lem_BCMP_Product_form}) and is summable. Further, for any $\omega \in \mathcal{A}$, choose a $\bar{T}_\delta$ such that for all $t \ge \bar{T}_\delta$,
\begin{equation}
 \frac{V_t(C_{b_\delta})}{t}(\omega) \le \sum_{s:n+r \ge b_\delta} \pi(s) + \frac{\delta}{4}.  
 \label{eq_bound}
\end{equation}

From \eqref{eq_epsilon}, we have the following,
\begin{align}
\epsilon_t & \leq \left| \sum_{s\in C_{b_\delta}}\left( \frac{R(s,t)}{t} - R^*(s) \right)\right| + \left| \sum_{s \in C_{b_\delta}^c}\left( \frac{R(s,t)}{t} - R^*(s) \right)\right| \nonumber \\
& \leq \left| \sum_{s\in C_{b_\delta}}\left( \frac{R(s,t)}{t} - R^*(s) \right)\right| + \sum_{s \in C_{b_\delta}^c} \left| \left( \frac{R(s,t)}{t} - R^*(s) \right)\right| \nonumber \\
& \leq \Lambda \left(\max_\phi f(\phi)\phi \right) \left|  \frac{V_t(C_{b_\delta})}{t} - \sum_{s\in C_{b_\delta}} \pi(s) \right| +  \nonumber \\
& \hspace{37mm} \sum_{s \in C_{b_\delta}^c} \left| \left( \frac{R(s,t)}{t} - R^*(s) \right)\right|. 
\label{eq_epsilon_reduced}
\end{align}
Note that the first equality follows from the triangle's inequality. The second step is again a result of application of triangle's inequality and the fact that $C_{b_\delta}^c$ is a finite set. The third step follows from \eqref{eqn_revenue_using_rrt} of Lemma \ref{lem_mr_derivation} and \eqref{eq_bound}. 

 Next, for any $\omega \in \mathcal{A}$, choose $\tilde{T}_\delta$ further large (if required), such that
 $$\text{ for all }s \in C_{b_\delta}^c, \ \left| \frac{R(s,t)}{t}(\omega) - R^*(s) \right| < \frac{\delta}{2 \lvert C_{b_\delta}^c\rvert}.$$ 
Thus, using \eqref{eq_bound} and the above argument, \eqref{eq_epsilon_reduced} can be further simplified to the following,
\begin{equation}
  \epsilon_t  \leq  \Lambda \left(\max_\phi f(\phi)\phi \right) \left( \frac{\delta}{2} + \frac{\delta}{2} \right) = \Lambda \left(\max_\phi f(\phi)\phi \right)\delta. \nonumber  
\end{equation}
Re-define $\tilde{\delta} := \nicefrac{\delta}{\Lambda (\max_\phi f(\phi)\phi )} $ and then we have the result. \eop 

\begin{lemma} \label{lem_monopoly}
Assuming \textbf{A.1} and \textbf{A.2}, for all $\phi \in [0, \phi_h]$, define $M(\phi) =  \frac{\Lambda}{2}  f(\phi)\phi $ and $m(\phi) =  \Lambda  f(\phi)\phi - e\phi$ then the function $M$ and $m$ is strictly concave function of $\phi$.

\end{lemma}

\noindent
\textbf{Proof of Lemma \ref{lem_monopoly}: }
Differentiating $M(\phi)$ with respect to $\phi$ we get,
$$M'(\phi) = \frac{\Lambda}{2} \left[ f(\phi) + f'(\phi) \phi\right]. $$ 
Under \textbf{A.1}, $f(\phi)$ is a strictly concave and decreasing function. From strict concavity of $f$, we have that $f'(\phi)$ is negative and decreasing. Thus, $f'(\phi)\phi$ is decreasing function of $\phi$. Hence, $M'(\phi)$ is  strictly decreasing function of $\phi$. Therefore $M(\phi)$ is a strictly concave function of $\phi$. Similarly, $m'(\phi) = M'(\phi) - e$ and hence, $m(\phi)$ is also a strictly concave function of $\phi$. \eop

\begin{lemma} \label{lem_monotonicity_of_h_1_function}
Under \textbf{A.1}, the function $d: \phi \mapsto 2 f(\phi) + \phi f'(\phi)$ is monotonically decreasing.
\end{lemma}

\noindent
\textbf{Proof of Lemma \ref{lem_monotonicity_of_h_1_function}: } 
Under \textbf{A.1}, $f(\phi)$ is a strictly concave and decreasing function. From strict concavity of $f$, we have that $f'(\phi)$ is negative and decreasing. Thus, $f'(\phi)\phi$ is decreasing function of $\phi$. Since, sum of two monotonically decreasing functions is again a monotonically decreasing function. Hence, we have the result.  \eop

\begin{lemma}
\label{lem_negative_partial_derivative}
The   following term in the numerator of the partial derivative 
$\frac{\partial h}{\partial \phi_i} $,
{\small $Z(\phi_i) := (f(\phi_i))^2+f(\phi_i)f(\phi_{-i})+f(\phi_{-i})f'(\phi_i)\phi_i$}  is decreasing w.r.t. $\phi_i$ under {\bf A}.1.
\end{lemma}

\noindent
\textbf{Proof of Lemma \ref{lem_negative_partial_derivative}:} From \textbf{A.1},  we know that $f$ is a decreasing function ($f'<0$) and hence, $f^2$ is also a decreasing function. Furthermore, $f$ is a strict concave function (so $f'$ is strictly decreasing) and hence $f'(\phi_i)\phi_i$ is also a decreasing function. Thus, we have the result.
\eop

We state the following lemma whose proof is immediate.
\begin{lemma}\label{lemma_oscillatory}
if $\phi \in \left(f^{-1}\left(\frac{2e}{\Lambda}\right], \phi_h \right)$ then $e\phi > \frac{\Lambda}{2}f(\phi)\phi > \Lambda f(\phi)\phi - e\phi$.
\eop
\end{lemma}

\end{document}

%% file: coop.tex
\section{Cooperation and Comparison}
\label{sec:comparison}

We now study the scenario in which the platforms seek to operate together. They combine their individual driver databases and attempt to serve the passengers together. The analysis of this scenario is exactly similar to that corresponding to the monopoly scenario (see Section~\ref{sec:monopoly}); by Theorem~\ref{thm_monopoly_opt}, the optimal price policy for limit combined system equals $\phi^* = \max\{\phifp, \phiunderbar\}.$  

\noindent{\bf Cooperation and monopoly:}
We immediately have an interesting observation about the comparison between monopoly (the platforms operate independently without interfering with each other) and cooperation (they operate together): the ratio $\nicefrac{e}{\Lambda}$ remains the same and thus {\it the optimal price (and hence the optimal matching revenue per platform) remains the same for both the systems.} This is because the value $\phifp$ depends solely upon the passenger response function $f$, while $\phiunderbar$ depends only upon the ratio
$\nicefrac{e}{\Lambda}$ and $f$. 

The above is the case with small $\beta$. In other words, when drivers are willing to wait for sufficiently long periods (as is usually the case with many practical scenarios), the usual economies of scale that we observe in queuing systems resulting from two or more system operating together is not seen with the double-sided queues.  
However the same is not the case when $\beta $ is sufficiently large. We consider an example in Figure~\ref{fig:monopoly_coop}, where we  numerically compute the optimal values of the MR using Lemma~\ref{lem_mr_derivation} and~\eqref{eqn_pi0_expression} for the cases with $\beta >0$. We now observe that the cooperation  makes a difference when $\beta$ is large. The platforms derive better optimal price and better matching revenue when they operate together. However the difference in optimal price across the two configurations is only about 4.6\% even for $\beta = 1$. On the other hand, the difference in matching revenues is more; it is around 17\% for $\beta = 0.5$. \textit{Thus one can conclude that cooperation does not improve significantly the revenue of the two platforms, unless the drivers are highly impatient.} 


\noindent{\bf Cooperation and competition:}
Next, we consider the competition between the two platforms. We evaluate and compare the equilibrium pricing and the corresponding matching revenues across the different passenger QoS models considered in Sections~\ref{sec:pi0} and~\ref{sec:overall_bp}, and also the monopoly (and the cooperation) setting. This comparison is performed for different choices of the price sensitivity function~$f.$

We begin with linear response function $f(\phi) = 1 - a \phi$ for some $a < 1/ \phi_h;$ this satisfies {\bf A}.1-2. After simple algebra, one can compute the following quantities of Theorems \ref{thm_monopoly_opt}-\ref{thm_epsilon_NE}
that define various equilibria as below: 
\begin{align*}
\phifp &= \frac{1 }{2a}
,  &\phido &= \frac{2  }{3a},
\\
 \phiunderbar &= \left( \frac{1}{a} - \frac{2e}{a\Lambda}\right),   &\phibar &=\left( \frac{1}{a} - \frac{e}{a\Lambda}\right), \\
\phimr &= \frac{1 - e/\Lambda}{2a} , 
&\eL &=  \frac{(\Lambda/e - e/\Lambda)^2}{4a}.
\end{align*}
We also consider a non-linear response function, $f(\phi)= 1 - (a \phi)^2$ for some $a < 1/ \phi_h$ that satisfies {\bf A}.1-2. For this function, the respective quantities are:
\begin{align*}
\phifp &= \frac{1 }{\sqrt{3}a}
,  &\phido &= \frac{1 }{\sqrt{2}a},
\\
 \phiunderbar &= \left(\frac{1}{a}\right) \sqrt{\left( 1 - \frac{2e}{\Lambda}\right)},   &\phibar &=\left(\frac{1}{a}\right) \sqrt{\left( 1 - \frac{e}{\Lambda}\right)}, \\
\phimr &= \left(\frac{1}{a}\right) \sqrt{\frac{1 - e/\Lambda}{3}} , 
&\eL &= \left(\frac{2}{3a}\right)\left(\frac{\Lambda}{e} - 1 \right) \sqrt{\frac{1 - e/\Lambda}{3}}.
\end{align*}
The immediate observation, which is strikingly visible from the above two sets of  expressions (which is also seen from Theorems \ref{thm_monopoly_opt}-\ref{thm_epsilon_NE}) is that the equilibrium prices as well as the normalized MR (MR at $\Lambda=1$) depend only upon  the arrivals-ratio $\nicefrac{e}{\Lambda}$ and response function $f$ and nothing else (the ratio is like the well known load-factor of the queuing systems). 

To derive further insights, 
we plot two case studies, one with linear response function is in Fig. \ref{fig:Lin_function} and the one with square response function is in Fig. \ref{fig:Square}. 
When  the arrivals-ratio $\nicefrac{e}{\Lambda}$ is small, the system is governed by effects of overwhelming number of passengers, and all the equilibria are highly sensitive to the above ratio. We observe this phenomenon for  both the response functions (left sub-figures of Fig. \ref{fig:Lin_function}  and  Fig. \ref{fig:Square}). When the arrivals-ratio is large, the equilibrium price  for   monopoly (and hence cooperation) as well as the duopoly-QoS-$\DA$ is insensitive to further increase in arrivals-ratio. 
However for the  duopoly driven by $\PB$, we observe the existence of equilibrium cycle (EC). At each such value of arrivals-ratio $\nicefrac{e}{\Lambda}$, we have a vertical line that represents the EC.  

We also capture the equilibrium MR for monopoly and duopoly (driven by $\PB$) scenarios in the right sub-figures. For the case with EC, we  plot the MR of a platform when it  prices at one  end point of the EC, while the opponent   offers  at the other end point (the end points would not matter as seen from proof of Theorem \ref{Thm_PB_EC}). 

We have several observations: a) interestingly the performance of the platforms is unaffected by the existence of the competition, for the cases with abundant passengers (or with $\nicefrac{e}{\Lambda}$ small);  b) the price and the  normalized  MR of monopoly (hence cooperation) equals that  corresponding to either of  the  two duopoly cases (up to   $\nicefrac{e}{\Lambda} \le 0.3 $ for linear case and up to $\nicefrac{e}{\Lambda} \le 0.45 $ for the square $f$, in right sub-figures); c) however the observation is very different when there is scarcity of passengers; cooperation significantly improves the MR of both the platforms for lager values of $\nicefrac{e}{\Lambda}$. 

In all, the price of anarchy is large when system operates with large arrival-ratios, while the same is negligible when it operates at medium or low arrival-rates.